%% file: main.tex
\documentclass{article}

\usepackage[english]{babel}

\usepackage[a4paper]{geometry}

\usepackage{amsmath}
\usepackage{graphicx}
\usepackage{xcolor}
\definecolor{foo}{HTML}{887B1B} 
\definecolor{Lichtgrijs}{HTML}{ADAFAF} 
\usepackage{url} 
\usepackage[colorlinks=true, citecolor=foo, linkcolor=foo, filecolor=foo, urlcolor=foo, pdfborder={0 0 0}]{hyperref}
\usepackage[square, numbers]{natbib}
\usepackage{booktabs}
\usepackage{pdfpages}
\usepackage{multirow}
\usepackage{makecell} 
\usepackage{tabularx}
\usepackage{colortbl}
\usepackage{pdflscape}
\usepackage{geometry}
\usepackage{caption}
\newcommand\edit[1]{{#1}}

\newcommand\unsure[1]{{#1}}
\newcolumntype{L}{>{\raggedright\arraybackslash}X}

\title{A Framework for Energy Management Modelling in Hubs for Circularity: Incorporating Industrial Symbiosis}

\usepackage{colortbl}      

\begin{document}
\begin{center}
{\huge A Framework for Sustainable Energy Management Modelling in Hubs for Circularity}
\end{center}

\renewcommand*{\thefootnote}{\fnsymbol{footnote}}
\vspace{-1mm}
\begin{center}
\mbox{\large Tobias Løvebakke Nielsen}\mbox{\large , Daniela Guericke}\footnote{Corresponding author}\mbox{\large, Alessio Trivella}\mbox{\large, Devrim Murat Yazan}\\[11pt]
{\small Industrial Engineering and Business Information Systems, University of Twente, Hallenweg 17, 7522 NH Enschede, The Netherlands. E-mails: \{t.l.nielsen, d.guericke, a.trivella, d.m.yazan\}@utwente.nl}\\[3pt]
\end{center}
\vspace{-1pt}
\renewcommand*{\thefootnote}{\arabic{footnote}}


\begin{center}
{\bf Abstract}
\end{center}

\begin{center}
\begin{minipage}[h]{0.90\columnwidth}

The concept Hubs for Circularity (H4C) represents integrated systems that combine efficient use of clean energy and circular economy to enhance resource efficiency within a region. H4C benefit from the geographical proximity of different industries within industrial zones and the surrounding urban and rural areas, allowing them to share resources, technology, and infrastructure. They reduce the use of virgin resources through Industrial Symbiosis (IS), where one company uses waste of another company as resource. 
Energy management is crucial in H4C, as these systems often integrate renewable energy sources, involve energy-intensive industries, include numerous energy consumers, and may rely on energy-based industrial symbiosis exchanges. 


This study presents a modelling framework for energy management in H4C, developed through a systematic literature review of related systems including energy-based IS. The framework provides a guideline for researchers and practitioners on which modelling aspects to consider when optimising the energy flows within a hub.   We argue that effective energy management in H4C requires combining conventional modelling aspects—such as objective functions, uncertainty, operational flexibility, and market participation—with IS-specific factors like the type of symbiosis, the degree of information sharing, and collaboration structures. In H4C with extensive IS, decentralised resource and energy exchanges often lead to similarly decentralised information flows and decision-making. Yet, our review shows a persistent reliance on centralised model structures, suggesting a path dependency rooted in traditional energy optimisation approaches. This highlights the need for models that better align with the distributed and collaborative nature of H4C systems.

\vspace{10pt}
\end{minipage}
\end{center}

\vspace{3pt}

\input{1_Introduction}

\input{2_Methodology}

\input{3_0_Results}

\input{3_1_waste_concept}

\input{3_3_obj_fun}

\input{3_4_Uncertainty}

\input{3_5_operational_response}

\input{3_6_energy_market}

\input{3_7_Centralised}

\input{3_2_Framework}

\input{4_Conclusion}

\section*{Acknowledgements}
This research is part of the project IS2H4C (Sustainable Circular Economy Transition: From Industrial Symbiosis to Hubs for Circularity) which has received funding from the European Union’s HORIZON Research and Innovation Actions programme under grant agreement number 101138473.
%
%

\section*{Declaration of generative AI and AI-assisted technologies in the writing process.}

During the preparation of this work the authors used ChatGPT for proofreading in order to improve the readability and language of the manuscript. After using this tool/service, the authors reviewed and edited the content as needed and take full responsibility for the content of the publication.

\input{AppendixTable}
{\small
\bibliographystyle{plainnat}
\bibliography{Bibliography}
}

\end{document}

%% file: 1_Introduction.tex
\section{Introduction}\label{sec:intro}
With its take-make-waste philosophy, the linear economy has been a main contributor to the ongoing climate, biodiversity, and pollution crises. Extraction and processing of material resources—including fossil fuels—account for 55\% of total greenhouse gas emissions, and if no action is taken, global resource extraction could rise by 60\% in 2060 compared to 2020 levels \cite{UNEP}. In response, both the UN and the EU recognise the circular economy as a key strategy to mitigate this development by keeping resources in use and minimising waste \cite{UNEP}. In a circular economy, systems are designed to minimise waste by using renewable energy and recirculating materials and waste, thereby reducing the demand for virgin finite resources \cite{EMF_butterfly}. 
Because waste generally holds low economic value, transporting it over long distances is often economically inefficient \cite{SHEU201465}. Instead, recirculation can be achieved through industrial symbiosis (IS), where waste from one production process is reused as input in another facilitated trough geographical proximity \cite{Energy_IS}. Through these symbiotic relationships, resources are being kept in the loops instead of leaking through \cite{EMF_butterfly}. 

\begin{figure}
    \centering
    \includegraphics[width=0.7\linewidth]{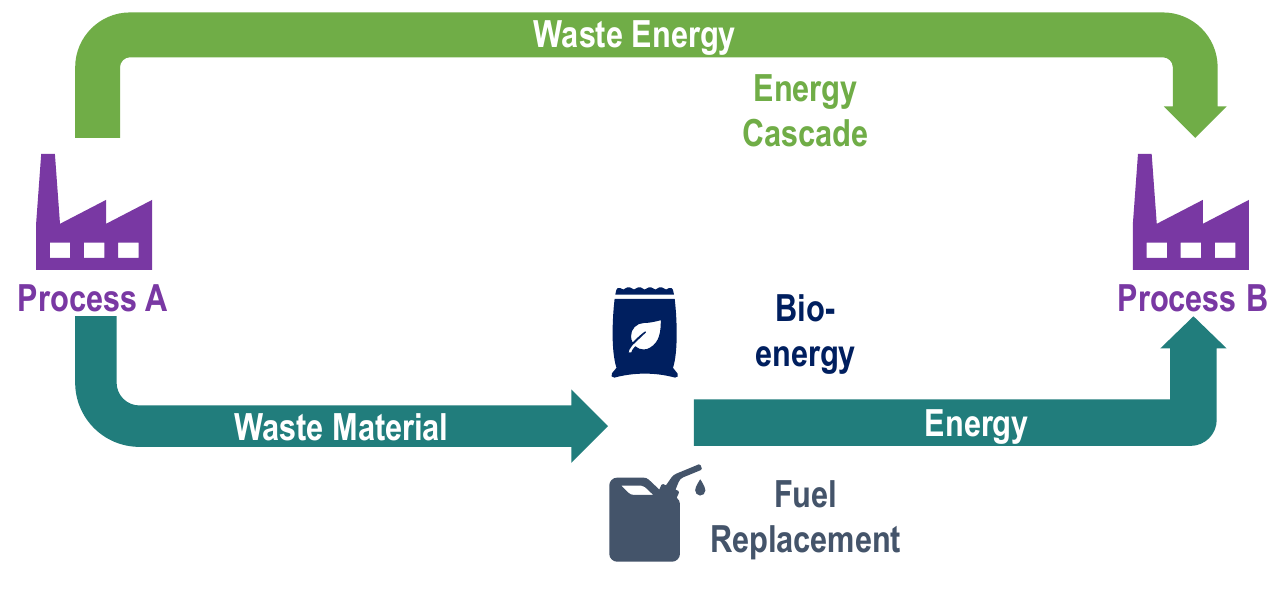}
    \caption{Type of energy-based IS as proposed by \citet{Energy_IS}}
    \label{fig:IS_energy}
\end{figure}

This idea is central to the Hubs for Circularity (H4C) concept, established by the Process4Planet partnership in 2021. H4C systematically apply IS in locally integrated industrial, urban and rural areas, aiming to create self-sustaining industrial ecosystems that slow, narrow, and close loops of energy, material, and data flows \cite{P4P_SRIA}. Classically, the waste being shared in an IS process is a solid material or resource. However, as \citet{Energy_IS} present in their work and as Figure \ref{fig:IS_energy} shows, IS can also be energy-based. They categorise three types of energy-based IS: energy cascade, bio-energy and fuel replacement. Energy cascade is when waste energy from one process is used as energy in another process, e.g. re-utilisation of waste heat. Bio-energy involves using bio-based waste as an energy source in a separate process, e.g\., conversion of manure into biogas. Lastly, fuel replacement covers waste material that replaces traditional fuel in fuel-based processes, e.g\., captured carbon-based fuel replacing natural gas in CHP plants. In addition to energy-based IS, energy supply and consumption play a vital role in H4C due to the presence of renewable energy sources, energy-intensive industries and other energy consumers in H4C. \unsure{This highlights the importance of effective energy management systems in H4C. In this study, we use the term energy management to refer to the modelling approach to optimise the IS and energy flows in order to obtain an energy-efficient system \cite{SHARMA2022120028}.}

The H4C energy system resembles conventional modern energy systems, such as energy communities, energy hubs, and positive energy districts, as it also leverages flexibility from processes, storage, and sector coupling to better utilise variable renewable energy. However, to fully enable H4C ecosystems, renewable energy management needs to be integrated with waste sharing to exploit synergies across sectors, leading to additional aspects that need to be considered in H4C energy management. In this paper, we address the question: \textit{Which aspects need to be considered when modelling energy management in H4C's?} Based on a systematic literature review on similar concepts including IS, we propose a framework for energy management in H4C. The framework serves a guideline for researchers and practitioners to identify relevant modelling aspects for energy management in H4C.

The remainder of the paper is structured as follows. Section \ref{sec:relatedwork} gives an overview of energy management reviews in related systems and outlines the contribution of our work in more detail. Section \ref{sec:methodology} presents the  methodology. Section \ref{sec:SLR} summarises and discusses the thematic review results. Based on the results, we propose the modelling framework in Section \ref{sec:framework}.  Section \ref{sec:conclusion} concludes the paper.

\section{Related Work and Contribution}\label{sec:relatedwork}

Various studies throughout the literature have characterised energy management in systems similar to H4C such as energy communities, eco-industrial parks, production systems and energy hubs. 

\textit{Energy communities} involve local residential actors pooling assets like rooftop PV, heat pumps, and batteries to cooperatively trade in electricity markets. In contrast, \textit{energy hubs} and \textit{positive energy districts} incorporate additional energy carriers such as heat and gas. Energy hubs serve as a general concept for all multi-carrier energy systems converting, storing and supplying various energy carriers \cite{MOHAMMADI20171512}, while positive energy districts focus on optimizing surplus electricity use within specific geographic areas \cite{TLEUKEN2025125324}. In these systems, energy demand typically defines the system boundary. However, in H4C, industrial symbiosis (IS) creates a more interconnected system, where processes not only consume but can also act as a waste provider of energy and non-energy resources for other processes.

A comparable approach is seen in \textit{eco-industrial parks}, which also promote sectoral collaboration in both energy and IS. However, IS relationships in eco-industrial parks are typically established through a top-down process \cite{EIP_WBG, HONG2020122853}. In contrast, H4C emphasize the development of spontaneous, bottom-up IS relationships that has historically led to more successful IS implementation \cite{GIBBS20071683}. Because such collaborations are initiated by the participating actors themselves, rather than being facilitated by a central entity, trust is more easily established. As a result, H4C systems are expected to feature a higher degree of IS and, consequently, a more interconnected operational structure than traditional eco-industrial parks.

Optimising energy flows in such complex and interconnected systems while addressing the inherent variability of renewable energy sources requires the development and application of advanced energy management models. How to model such integrated systems is already an active research field, and several review articles have sought to synthesise the modelling trends. In the context of \textit{production environments}, \citet{BANSCH2021107456} conducted a systematic literature review (SLR) that maps out the types of energy management problems studied and the solution approaches adopted. Similarly, \citet{Xavier_SCM_EM} present a modelling-focused SLR within the field of supply chain management, where energy is primarily considered as an external upstream supply. Based on this framing, they also propose a conceptual framework to integrate energy considerations into supply chain models.

 Several studies have focused on \textit{energy hub} concepts. \citet{MOHAMMADI20171512} present a conceptual energy hub model and review the energy carriers and technologies covered in the literature. \citet{NOZARI2022103972} offer a review focused on energy storage, also proposing a model structure and summarising the storage technologies and solution approaches used. A more comprehensive modelling perspective is presented in \citet{SHABANPOURHAGHIGHI20226164} and \citet{NOZARIAN2024103834}. \citet{SHABANPOURHAGHIGHI20226164} provide a thematic analysis of the applied objective functions, model approaches, demand response approaches and uncertainty approaches, while \citet{NOZARIAN2024103834} include a thematic analysis of energy technologies, general model structure and KPI's of building cluster meso-scale energy hubs.  Lastly, \citet{KIANIMOGHADDAM2023128263} examine how uncertainty and solution methods are handled across various energy hub studies.

However, IS is notably absent from these reviews \cite{BANSCH2021107456, Xavier_SCM_EM, MOHAMMADI20171512, NOZARI2022103972, SHABANPOURHAGHIGHI20226164, NOZARIAN2024103834, KIANIMOGHADDAM2023128263}. In contrast, Fraccascia et al. \cite{Energy_IS} focus explicitly on energy-based IS, categorising different types and identifying associated drivers, barriers, and enablers. Turken and Geda \cite{TURKEN2020104974} examine the synergies between IS and supply chains, including energy flows. Yazici Alakaş and Eren \cite{yazici_OR_IS} provide a comprehensive review of operational research models applied to IS decision-making, although most of these do not include energy systems. Butturi et al. \cite{BUTTURI2019113825} explore the potential of integrating renewable energy sources into eco-industrial parks, proposing a framework that highlights the key elements addressed in eco-industrial park studies. More recently, \citet{RAMIRDTCERTEZA2025115377} conducted a detailed review on the integration of waste heat into energy-based IS networks, covering aspects such as objective functions, uncertainty, spatial distance, and process types.

Despite this extensive existing literature, a research gap remains. To the best of our knowledge, no comprehensive review has yet addressed energy management modelling in sector-coupled energy systems that integrate IS principles, as seen in emerging concepts such as H4C. Furthermore, to the best of our knowledge, no previous work has systematically categorised the operational IS decisions embedded in optimisation models of multi-carrier energy concepts according to how circular they are. This paper aims to address these gaps by:
\begin{enumerate}
    \item Proposing a comprehensive modelling framework that outlines the aspects that should be considered in energy management models for H4C systems, including the specific implications of industrial symbiosis on the applied modelling approach;
    \item Conducting a systematic literature review of energy management models relevant to H4C, including thematic reviews of IS types, objective functions, uncertainty factors, energy market participation, operational flexibility, and overall modelling approaches.
    \item Categorising the operational decisions of multi-carrier energy system optimisations according to the waste hierarchy levels, as defined in \citet{KIRCHHERR2017221}.
\end{enumerate}

By doing so, this study contributes to a clearer understanding of energy management modelling in H4C contexts and offers structured guidance for both practitioners and researchers engaged in developing such models.

%% file: 2_Methodology.tex
\section{Methodology} \label{sec:methodology}
The purpose of the review is to answer the following research question: 
\textbf{\textit{Which aspects need to be considered when modelling energy management in H4C's?}}

\begin{figure}[ht]
    \centering    \includegraphics[width=0.6\textwidth]{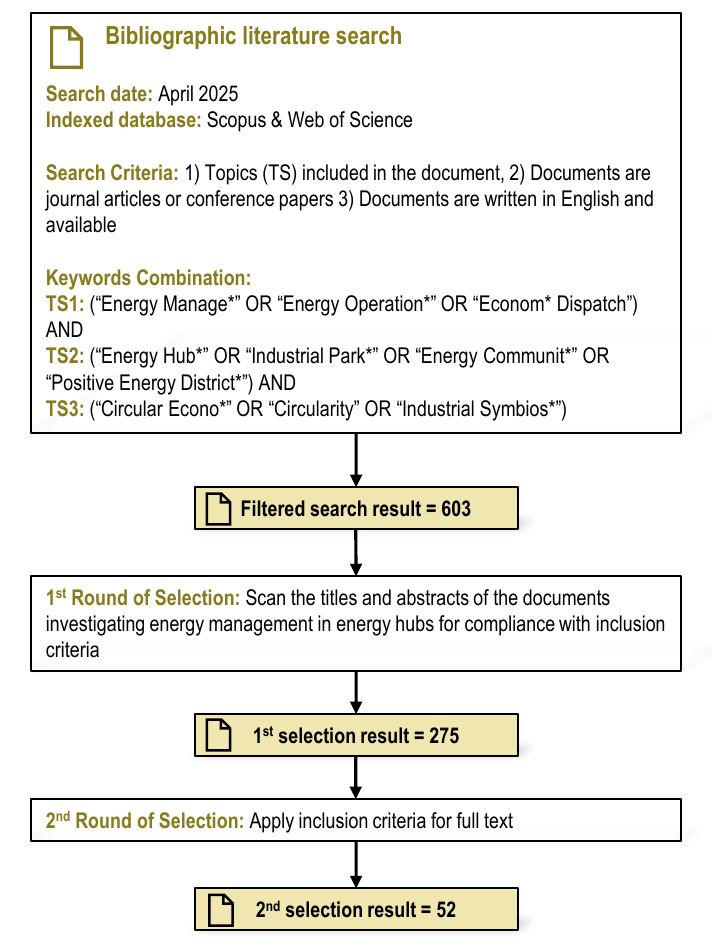}
    \caption{SLR Methodology}
    \label{fig:SLR_methodology}
\end{figure}

 We conducted a systematic literature review (SLR) to gain insights on requirements and modelling aspects for energy management that are relevant for H4C. The SLR is based on the item checklist from the Prisma 2020 statement \cite{Pagen160}. Figure \ref{fig:SLR_methodology} shows the general methodology.  On the basis of the review, we propose a modelling framework to highlight which modelling aspects are important to answer the research question and to identify research gaps.

We conducted the search in April 2025 using the academic bibliometric databases Scopus and Web of Science. We designed the query around three main search criteria and applied two selection rounds. The first search criterion was the search query. It combined two primary thematic areas: energy management and hubs for circularity (H4C), connected by an “AND” operator to ensure that only studies addressing both themes were included. To represent energy management, the search included the terms “energy operations" and “economic dispatch" (TS1). The second theme was more challenging to define, as the term “H4C" was introduced only in 2021 by the Process4Planet partnership \cite{P4P_SRIA}, and to our knowledge, no specific literature on energy management within H4C currently exists at the time of submission for this study. Therefore, we defined two additional keyword sets (TS2 and TS3) representing related concepts to H4C that also incorporate circular economy concepts. TS2 includes \textit{energy hub, industrial park, energy communities} and \textit{positive energy districts}. TS3 cover circular economy-related terms \textit{circular economy, circularity} and \textit{industrial symbiosis}. In addition to the search query, we defined the second criterion which is that the considered documents have to be either peer-reviewed journal articles or conference papers. The third criterion was that the document had to be written in English and be available for assessment. The review was conducted by the first author independently.

\begin{figure}
    \centering
    \includegraphics[width=0.45\textwidth]{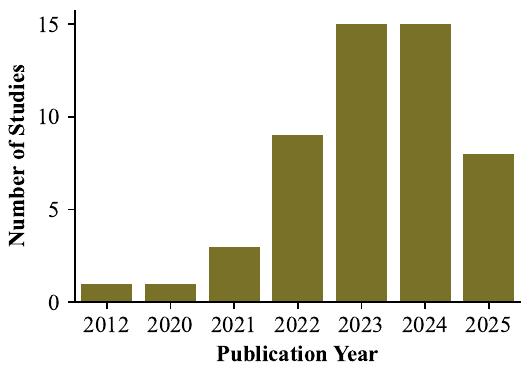}
    \caption{Study distribution with regards to publication year}
    \label{fig:dist_year}
\end{figure}

The initial search yielded 603 documents. In the first selection round, we screened titles and abstracts based on two inclusion criteria:
\begin{enumerate}
    \item The considered document has to implement a mathematical optimisation problem
    \item The mathematical optimisation problem has to consider multiple energy carriers
\end{enumerate}
These criteria were applied to ensure the selected studies address operational-level decision-making and exclude studies focused solely on organizational or governance aspects, as discussed in \cite{ABDELAZIZ2011150}. Furthermore, given that H4Cs typically involve integration across several energy carriers, this focus ensures alignment with the modelling needs of H4C systems. The first screening reduced the search to 275 studies. In the second round, we screened the entire document for compliance with the inclusion criteria. After this round, 52 studies remained for in-depth review. The distribution of these studies over time is shown in Figure \ref{fig:dist_year}, indicating that this is a relatively new area of research, with increasing interest in recent years. Notably, only one study was published before 2020. 

%% file: 3_0_Results.tex
\section{Review of Energy Management Aspects Relevant for H4C} \label{sec:SLR}

\edit{This section presents the results of the SLR on how various energy management modelling aspects have been applied in concepts similar to H4Cs. In Section \ref{sec:waste}, we start by giving an overview of application cases of energy-based IS in relation to the considered context (energy hubs, industrial park and energy community). In addition, we categorise the optimisation decisions based on the circular economy waste hierarchy \citep{KIRCHHERR2017221}. Section \ref{sec:SLR_obj} and Section \ref{sec:uncertainty} review the objectives and uncertainty factors addressed in the literature. Section \ref{sec:flex} covers flexibility measures, including operational response (Section \ref{sec:DR}) and storage (Section \ref{sec:storage}). Section \ref{sec:market} discusses how different markets are integrated into the models. Sections \ref{sec:central} and \ref{sec:decentral} classify studies by centralised or decentralised modelling approaches and discuss the applied methodologies. To broaden the decentralised modelling perspective, eight relevant studies from related fields are also included. Table \ref{tab:app} in Appendix \ref{App:tab} provides a general overview of all reviewed references including their considered modelling aspects.}

%% file: 3_1_waste_concept.tex
\subsection{Energy-based IS and Categorisation on Waste Hierarchy} \label{sec:waste}

As mentioned in Section \ref{sec:intro}, \edit{IS itself can be energy-based in terms of energy cascades, bio-energy and fuel replacement. These three IS-types are present in most of the reviewed literature. However, some types increase the circularity of systems more than others. This is the essence of the 9R circularity hierarchy, where the lower R a circularity strategy has, the more it increases the circularity of a system \cite{KIRCHHERR2017221}. In Table \ref{tab:WH_IS}, we categorise the considered optimisation decisions of the reviewed studies according to the circularity hierarchy levels together with the applied energy-based IS type. In this study, we only consider a paper to apply IS if its optimisation decisions explicitly incorporate waste flows.} For example, producing hydrogen via electrolysis may be preferable from a circularity design perspective compared to producing hydrogen from natural gas. However, since hydrogen is intentionally produced, it is not a waste product and thereby not considered as operational IS in this study. Additionally, Table \ref{tab:WH_IS} also categorises the studies into the three main academic concepts that are applied throughout the considered literature: \textit{energy hubs, industrial parks} and \textit{energy communities}. Despite being part of query, no \textit{positive energy district}-study fulfilled the inclusion criteria.

Table \ref{tab:WH_IS} shows that all 52 reviewed studies apply the reduce strategy (R2) by optimising resource usage and thereby narrowing resource loops within their defined system boundaries \cite{Bocken03072016}. However, only half of these studies (26) include IS optimisation decisions, where waste is re-utilised in downstream processes. The IS flow decisions lead to a more interconnected supply, and given the balancing requirements of energy systems, increases the complexity of energy management. Nevertheless, we still included the remaining 26 studies without optimisation of IS decisions (marked with X in Table \ref{tab:WH_IS}), as they still offer valuable insights into the various modelling aspects of energy management in interconnected energy systems.

\setlength{\tabcolsep}{4pt}
\begin{table*}[h!]\centering 
\footnotesize
\rowcolors{2}{gray!20}{white} 
\begin{tabularx}{\textwidth}{c| c c c c c c c c c c p{0.2cm} L}
\rowcolor{white}
\multicolumn{1}{c}{\textcolor{foo}{\textbf{Concept}}} &
\multicolumn{10}{c}{\textcolor{foo}{\textbf{Waste Hierarchy Levels}}} &
 & 
\multicolumn{1}{c}{\textcolor{foo}{\textbf{Papers}}} \\
\cmidrule{1-11} \cmidrule{13-13}
&\makecell[t]{\scriptsize \textbf{R0} \\ \tiny Re-\\\tiny fuse} & \makecell[t]{\scriptsize \textbf{R1} \\ \tiny Re-\\\tiny think} & \makecell[t]{\scriptsize \textbf{R2} \\ \tiny Re-\\\tiny duce} & \makecell[t]{\scriptsize \textbf{R3} \\ \tiny Re-\\\tiny use} & \makecell[t]{\scriptsize \textbf{R4} \\ \tiny Re-\\\tiny pair} & \makecell[t]{\scriptsize \textbf{R5} \\ \tiny Re-\\\tiny furbish} & \makecell[t]{\scriptsize \textbf{R6} \\ \tiny Re-\\\tiny manu- \\\tiny facture} & \makecell[t]{\scriptsize \textbf{R7} \\ \tiny Re-\\\tiny purpose} & \makecell[t]{\scriptsize \textbf{R8} \\ \tiny Re-\\\tiny cycle} & \makecell[t]{\scriptsize \textbf{R9} \\ \tiny Re-\\\tiny cover} &&\makecell[c]{52} \\ \cmidrule{1-13}
\cellcolor{white}\textcolor{foo}{\textbf{Energy Hub}} &   &    & EC & EC &    &    &    &    &    &    && \cite{ZHANG2023127643}, \cite{REN2023119499},  \cite{CHINESE2022124785}, \cite{CAI2025136178}, \cite{PU_10637267}, \cite{WEN2025134696} \\
\cellcolor{white} &   &    & FR  &    &    &    &    &    &  FR &   && \cite{WANG2024129868}, \cite{LI2024119946}, \cite{KHALIGH2023117354}, \cite{ZUO2025135413}\\
\cellcolor{white}&   &    & BE &    &    &    &    &    &  BE &     && \cite{GOH2023117484}, \cite{fan10878712} \\
\cellcolor{white}&   &    & BE  &   &    &    &    &    &  BE & BE  && \cite{ZHAO2023101340}, \cite{LIU2024133528}\\
\cellcolor{white} &   &    & EC,FR & EC &    &    &    &      & FR&    && \cite{LI2022105261}\\
 \cellcolor{white}&   &    & EC,BE & EC &    &    &    &      & BE&    && \cite{Zhang_DTU}\\
\cellcolor{white}&   &    & FR,BE  &    &    &    &    &    &  FR,BE & BE  && \cite{MOSTAFAVISANI2023116965} \\
\cellcolor{white}&   &    & X  &    &    &    &    &    &  &  && \cite{WANG2023121141}, \cite{MUGHEES2023121150}, \cite{WANG2023126288}, \cite{WANG20234631}, \cite{GHASEMIMARZBALI2023120351}, \cite{SOHRABITABAR2022134898}, \cite{DAVOUDI2022107889}, \cite{Dwijendra}, \cite{ALABI2022118997}, \cite{LIN2021117741}, \cite{XU2024123525}, \cite{CHEN2024120734}, \cite{su16125016}, \cite{LI2024129784}, \cite{MOHTAVIPOUR2024132880}, \cite{REN2025135422}, \cite{SALYANI2025106245} \\
\cmidrule{1-13} 
\cellcolor{white}\textcolor{foo}{\textbf{Industrial}}&   &    & EC & EC &    &    &    &    &    &    && \cite{GALVANCARA2022117734}\\
\cellcolor{white}\textcolor{foo}{\textbf{Park}} &   &    & FR  &    &    &    &    &    &  FR &    && \cite{QIAN2024110504}, \cite{HOU2024130617}, \cite{chen_iet_nash} \\
\cellcolor{white}&   &    & BE &    &    &    &    &    &    & BE && \cite{MENEGHETTI2012263} \\
\cellcolor{white}&   &    & EC,BE  & EC  &    &    &    &    &  & BE  && \cite{Kantor_Frontier}\\
\cellcolor{white}&   &    & FR,BE  &    &    &    &    &    &  FR,BE & BE  && \cite{Chen_Wei} \\
\cellcolor{white}&  &  & X   &    &    &    &    &    &    &  && \cite{en16248041}, \cite{PANG2023125201}, \cite{GUO202133039}, \cite{TOSTADOVELIZ2024123389}, \cite{WEI2022121732}, \cite{HWANGBO2022122006}, \cite{DONG20251122} \\
\cmidrule{1-13}
\cellcolor{white}\textcolor{foo}{\textbf{Energy}}&  &  & BE &    &    &    &    &    & BE &  && \cite{pedro10161834} \\
\cellcolor{white}\textcolor{foo}{\textbf{Community}} &  &    & BE &    &    &    &    &    &    & BE && \cite{MALDET2024100106} \\
\cellcolor{white}&  &  & X  &    &    &    &    &    &    &  && \cite{LIU2024118204}, \cite{ZHENG2025101098} \\
\bottomrule
\end{tabularx}
\caption{Table of applied energy based industrial symbiosis and at what waste hierarchy level. Energy hubs covers here a broad range of concept incorporating multiple energy carriers, e.g. multi-carrier energy networks and micro-grids.  \textbf{EC}: Energy Cascade, \textbf{FR}: Fuel Replacement, \textbf{BE}: Bio-Energy, \textbf{X}: No IS optimisation decisions}
\label{tab:WH_IS}
\end{table*}
\textbf{Energy hubs.}
Among the reviewed studies, energy hubs are the most applied concept, appearing in 34 out of 52 studies. \edit{Among these, \citet{CAI2025136178, PU_10637267, WEN2025134696, CHINESE2022124785,REN2023119499, ZHANG2023127643} all apply energy cascading by reusing waste heat (R3) to supply a heat demand. \citet{CAI2025136178}, \citet{PU_10637267} and \citet{WEN2025134696} recover waste heat from electrolysers, a power plant and from carbon capture, respectively. \citet{CHINESE2022124785} reuse waste energy from the chillers of a cheese factory, while \citet{ZHANG2023127643} and \citet{REN2023119499} recover waste heat from fuel production to supply combined cooling, heat and power (CCHP) loads. Additionally, \citet{ZHANG2023127643} reuse waste water from fuel cell operation as input to an electrolysis process.}

\citet{WANG2024129868},  \citet{LI2024119946}, \citet{ZUO2025135413} and \citet{KHALIGH2023117354} also apply circular economy strategies, but in their studies, they recycle captured carbon to produce fuels, thus applying fuel replacement (R8). Both in \citet{WANG2024129868} and \citet{LI2024119946}, the captured carbon is converted into gas to both supply a combined heat and power (CHP) load and a gas load. In addition, \citet{LI2024119946} also supply a cooling load. Meanwhile, both \citet{ZUO2025135413} and \citet{KHALIGH2023117354} optimise a system where captured carbon is converted into fuels to be sold to external markets. Besides fuel selling, \citet{ZUO2025135413} optimise energy flows related to heat, power and gas loads, while \citet{KHALIGH2023117354} optimise energy flows of power, hydrogen and gas.

Energy-based IS in terms of bio-energy is investigated by \citet{GOH2023117484}, \citet{ZHAO2023101340}, \citet{LIU2024133528} and \citet{fan10878712}. In \citet{GOH2023117484}, various processes in an eco-industrial park generate biogenic waste, from which hydrogen can be extracted. Combined with solar power, the extracted hydrogen is then supplied back to the energy-consuming plants, creating a closed-loop energy system (R8). \edit{\citet{fan10878712} apply a similar IS strategy by converting manure into biogas (R8) used to supply a CHP load.} On the other hand, Zhao et al. \cite{ZHAO2023101340} optimise a system processing rural and urban waste partly into biogas for gas supply (R8). The remaining waste part is incinerated to generate heat and power (R9). \edit{Similarly, \citet{LIU2024133528} apply incineration of waste through pyrolysis of solid biomass waste (R9) and converts lifestock waste into biogas through anaerobic digestion  to supply heat, power and gas (R8). }

\citet{LI2022105261}, \citet{Zhang_DTU}, and \citet{MOSTAFAVISANI2023116965} apply combined IS strategies. In the study by \citet{LI2022105261}, a gas-fired power plant produces methane through methanation of captured carbon (R8), while recovering waste heat from methanation to meet a CHP demand (R3). Similarly, \citet{Zhang_DTU} utilize waste heat from gas production (R3 \& R8) within a biomass-based energy hub. On the other hand, Sani et al. \cite{MOSTAFAVISANI2023116965} optimize a system that both incinerate municipal solid waste (R9) and converts it into syn- and biogas (R8) to supply heat, power and hydrogen demands.

The remaining 17 studies optimise sector-coupled systems without IS flow decisions. Several focus on CCHP systems, including \citet{ALABI2022118997}, \citet{LIN2021117741}, \citet{su16125016}, \citet{LI2024129784} and \cite{REN2025135422}, while \citet{SOHRABITABAR2022134898} and \citet{DAVOUDI2022107889} also include gas supply. Gas is likewise considered in microgrids with electric vehicles in \cite{WANG2023121141} and in combination with a CHP supply in \cite{Dwijendra}. Water and hydrogen are supplied alongside CHP in \cite{GHASEMIMARZBALI2023120351} and \cite{XU2024123525}, respectively. \citet{WANG2023126288}, \citet{WANG20234631} and \citet{SALYANI2025106245} examine power–gas and power–hydrogen coupling, while \citet{CHEN2024120734} focus on solar thermal and cooling. Lastly, \citet{MUGHEES2023121150} and \citet{MOHTAVIPOUR2024132880} optimise a CHP supply in an industrial hub and microgrid, respectively.

\textbf{Industrial parks.}
The industrial park studies show a similar pattern, where also around half of them consider operational IS decisions (seven out of 14). \citet{GALVANCARA2022117734} consider an eco-industrial park, where multiple plants recover waste heat and thereby apply energy cascading (R3). \citet{QIAN2024110504}, \citet{HOU2024130617} and \citet{chen_iet_nash} all consider fuel replacement similarly to the energy hub case by converting captured carbon into gas (R8), while also producing hydrogen. \citet{QIAN2024110504} optimise a CHP supply, while both \citet{HOU2024130617} and \citet{chen_iet_nash} optimise a CCHP supply. \citet{MENEGHETTI2012263}, \citet{Kantor_Frontier} and \citet{Chen_Wei} all consider eco-industrial parks in which biomass is incinerated to produce energy (R9). In addition to incineration, \citet{Kantor_Frontier} also reuse waste heat in their study (R3), while \citet{Chen_Wei} uses biogas as a replacement of traditional fuels (R8).

The remaining seven industrial park studies optimise sector-coupled systems without operational IS decisions. Six focus on CHP supply, applied in contexts such as paper-making \cite{en16248041}, carbon-neutral operations \cite{WEI2022121732}, combined hydrogen gas production \cite{GUO202133039}, and petrochemical facilities \cite{HWANGBO2022122006}. \citet{DONG20251122} also optimise the supply of cooling besides CHP supply, while \citet{PANG2023125201} also optimise the supply of gas and cooling. Lastly, \citet{TOSTADOVELIZ2024123389} consider a combined power and hydrogen supply.

\textbf{Energy communities.}
The last considered concept throughout the studies is energy communities. \citet{pedro10161834} consider the conversion of agricultural waste into biogas used for supplying a CHP load (R8), while \citet{MALDET2024100106} investigate waste incineration for a CCHP supply in local sustainable communities (R9). Lastly, both \citet{LIU2024118204} and \citet{ZHENG2025101098} do not implement operational IS decisions. \citet{LIU2024118204} optimise a CCHP supply for a rural community, while \citet{ZHENG2025101098} optimise a CHP supply. 

\vspace{5pt}
All studies for energy hubs, industrial parks and energy communities included in the review optimise multiple interconnected energy carriers to achieve a more efficient energy supply and consumption. By including the IS synergies such as the ones presented in the 26 out of 52 studies, the optimisation goes one step further in considering also waste flows between companies. This results in a more holistic view that optimises energy carriers including parts of the surrounding ecosystem. An effective energy management in H4C needs to consider these waste flows in addition to the traditional energy management perspective.  

Furthermore, we can conclude that the different types of energy-based IS happen on different levels of the waste hierarchy, apart from R2 since all studies aim at reducing consumption in the first place. Energy cascades are related only to level R3 - Reuse. Fuel replacement is related to level R8 (Recycle), primarily through the re-utilisation of captured carbon. Bio-energy IS spans levels R8 or R9 (Recover), either through the incineration of bio-waste or the conversion of waste into biogas.

%% file: 3_3_obj_fun.tex
\subsection{Objective Function} \label{sec:SLR_obj}

Depending on the factors included in the objective function, an energy management method will lead to different results. Thus, the objective function governs how resources and energy are being shared. This section provides an overview of the objective functions in the reviewed literature. 

\begin{table*}[h!]\centering
\footnotesize
\rowcolors{2}{white}{gray!20} 
\begin{tabularx}{\textwidth}{c c c c c c c  X}
\multicolumn{7}{c}{\textcolor{foo}{\textbf{Objective Function Elements}}} & \multicolumn{1}{c}{\textcolor{foo}{\textbf{Papers}}} \\
\cmidrule(r){1-7} \cmidrule(r){8-8}
\rowcolor{white}
Economic & Emissions & Curtailment  & Degradation & Grid Import& Energy & DRU & \multicolumn{1}{c}{52}\\ \midrule
X &&&&&&& \cite{LIU2024118204}, \cite{MALDET2024100106}, \cite{WANG2023121141}, \cite{WANG2023126288}, \cite{GHASEMIMARZBALI2023120351}, \cite{SOHRABITABAR2022134898}, \cite{TOSTADOVELIZ2024123389}, \cite{HWANGBO2022122006}, \cite{Dwijendra}, \cite{GUO202133039}, \cite{MENEGHETTI2012263}, \cite{WANG2024129868}, \cite{WEI2022121732}, \cite{MOHTAVIPOUR2024132880}, \cite{LIU2024133528}, \cite{ZUO2025135413}, \cite{pedro10161834} \\
X &X  &&&&&&\cite{PANG2023125201}, \cite{LIN2021117741}, \cite{Chen_Wei},  \cite{ZHANG2023127643}, \cite{KHALIGH2023117354}, \cite{HOU2024130617}, \cite{ZHAO2023101340}, \cite{en16248041}, \cite{MOSTAFAVISANI2023116965}, \cite{CHINESE2022124785}, \cite{GALVANCARA2022117734},  \cite{XU2024123525}, \cite{CHEN2024120734}, \cite{Kantor_Frontier}, \cite{CAI2025136178}, \cite{PU_10637267}, \cite{DONG20251122}, \cite{WEN2025134696}, \cite{ZHENG2025101098}, \cite{SALYANI2025106245}, \cite{fan10878712}\\
X & X & X &&&&& \cite{QIAN2024110504}, \cite{LI2024119946} \\
X & X & X & X&&&& \cite{ALABI2022118997} \\
X && X & X &&&& \cite{WANG20234631} \\
X &&& X &&&& \cite{MUGHEES2023121150}, \cite{Zhang_DTU} \\
X && X && X &&& \cite{GOH2023117484} \\
X &&& && X && \cite{DAVOUDI2022107889}\\
X & X &&&& X && \cite{REN2023119499}, \cite{su16125016}, \cite{LI2024129784}, \cite{REN2025135422} \\
X&X& X& && & X& \cite{chen_iet_nash}\\
&&& && X  && \cite{LI2022105261}\\
\bottomrule
\end{tabularx}
\caption{Objective functions considered in the reviewed literature. DRU: Demand Reponse Utility}
\label{tab:obj}
\end{table*}

Table \ref{tab:obj} shows that 51 out of 52 studies include an \textit{economic} metric in their objective function. Of these, 17 studies optimise their decision variables solely based on economic metrics. The objective function in \citet{LIU2024118204}, \citet{Dwijendra} and \citet{WANG2024129868} consist of the operational costs. \citet{WANG2023121141} extends this by also considering electric vehicle destruction cost. Other studies include also trading costs alongside operational costs \cite{GHASEMIMARZBALI2023120351,SOHRABITABAR2022134898,TOSTADOVELIZ2024123389,WANG2023126288,GUO202133039,MOHTAVIPOUR2024132880, LIU2024133528, ZUO2025135413}. Additionally, some studies include resource disposal costs \cite{MALDET2024100106} while others introduce investment costs \cite{HWANGBO2022122006, WEI2022121732, MENEGHETTI2012263, pedro10161834} in their objective functions. 

The remaining 34 studies that include economic aspects also consider non-economic metrics. Of these, 21 \textit{combine economic and emissions-related objectives}. \citet{LIN2021117741}, \citet{PANG2023125201}, \citet{ZHENG2025101098}, \citet{Chen_Wei} and \citet{ZHANG2023127643}, \citet{PU_10637267}, \cite{SALYANI2025106245} utilize a Pareto front with an economic and an emission objective dimension to highlight trade-offs between minimizing costs and reducing emissions. As the economic dimension, \citet{LIN2021117741}, \citet{PANG2023125201} and \citet{ZHENG2025101098} include operating and investment costs, while \citet{Chen_Wei}, \citet{ZHANG2023127643}, \citet{PU_10637267}, \cite{DONG20251122} and \citet{SALYANI2025106245} include operating and trading costs. As a third dimension, \citet{SALYANI2025106245} also include the risk measure conditional value at risk, while \citet{DONG20251122} include hydrogen price as a third dimension. 
In the remaining 14 articles that combine economic and emissions metrics, the emission metric is internalised as a penalty cost element. These are considered either alongside operational and trading cost \cite{KHALIGH2023117354,HOU2024130617,ZHAO2023101340,en16248041,CHEN2024120734, CAI2025136178, PU_10637267} or operational and investment cost \cite{MOSTAFAVISANI2023116965,CHINESE2022124785,GALVANCARA2022117734,XU2024123525,Kantor_Frontier, WEN2025134696, fan10878712}.

Apart from emissions, non-econometric aspects considered are \textit{curtailment, degradation of batteries and shares of certain energy types}.  \citet{QIAN2024110504} combine curtailment and emissions into a cost objective also covering operating and trading costs. Similarly, curtailment is also part of the cost objective together with operating, trading and investment costs in the study by \citet{LI2024119946}. Additionally, they include an emission aspect as a second Pareto Front dimension. \citet{ALABI2022118997} apply a similar objective to \cite{QIAN2024110504}, but they also include the battery {degradation} cost associated with providing operational flexibility. \citet{WANG20234631} also internalise degradation of batteries as a cost together with operating, trading and curtailment cost. However, they do not include emissions in their objective. \citet{MUGHEES2023121150} and  \citet{Zhang_DTU} consider cost for degradation but not curtailment.  \citet{MUGHEES2023121150} also consider investment costs. 
\citet{GOH2023117484} have a four-dimensional objective function featuring annualized operating costs, unused energy, and shares of imported electricity and hydrogen.

The studies \cite{DAVOUDI2022107889, REN2025135422, LI2024129784, REN2023119499, su16125016, chen_iet_nash, LI2022105261} all included metrics related to the \textit{amount of used energy}.  \citet{DAVOUDI2022107889} use a multi-objective setting with a three parts combining investment cost, operating and trading cost and a metric for energy reliability with regards to failures. \citet{REN2025135422} has similar economic and emissions dimensions, but they include the energy deficit of their system as their energy dimension.  \citet{LI2024129784}, \citet{REN2023119499}, \citet{su16125016} and \citet{chen_iet_nash} all consider economic, emissions and energy related terms in their objective. \citet{LI2024129784} construct a three dimensional Pareto front of operating and trading cost, carbon emissions and the amount of energy being used. Similarly,  \citet{REN2023119499} implement a multi-objective setting looking at operating costs, emissions and conversion efficiency of a solar thermochemical plant. \citet{su16125016} on the other hand consider a single weighted objective of economic, emissions and energy savings compared to a reference solution. \citet{chen_iet_nash} consider a weighted objective consisting of operating and trading costs, carbon trading costs, curtailment and utility of providing demand response (DRU). 

\citet{LI2022105261} are the only ones that do not consider economic metrics in the reviewed literature. They solely consider the amount of gas supplied to a power plant in their objective. 

\edit{Most of the studies optimise a system with a joint objective. Only \citet{ZHENG2025101098} consider the allocation of shared benefits among the included stakeholders. In their case, they distributed the benefits from providing demand response based on an overall weighted average across the entire production period. However, due to the fluctuating nature of renewable energy generation throughout a day, a simple weighted average may not accurately reflect each participant’s time-dependent contribution to the shared benefit. This issue is critical not only for the distribution of profits and costs but also for ensuring fair allocation of renewable energy. As corporate interest in sustainability and circularity continues to grow \cite{wef2025circular}, stakeholders in a hub might also compete on utilization of the available renewable energy to fulfil their targets. Therefore, the question of how to fairly allocate renewable energy to consumers actively contributing to its uptake within systems similar to H4C remains an open research gap.}

In conclusion, economic performance is the main objective within energy management of systems similar to H4C, with 51 out of 52 studies including economic metrics in their objective. However, since H4Cs focus on the circularity of a system \cite{P4P_SRIA}, economic efficiency is not the sole optimised focus. This is also reflected throughout the reviewed literature on similar concepts as H4C, where a majority of the studies include additional objective dimensions alongside economic metrics. Among these, emissions reduction was the most common objective term, included in 29 of the studies. In contrast, resource efficiency metrics, such as curtailment and degradation, receive comparatively less attention, appearing in only six and four studies, respectively. 

The inclusion of multiple metrics in multi-objective setting illustrates how different optimization preferences influence system behaviour and resource-sharing dynamics within hubs. This highlights the importance of considering objective function design when modelling energy management strategies for H4C.

%% file: 3_4_Uncertainty.tex
\subsection{Uncertainty Factors} \label{sec:uncertainty}
The energy management in H4C is influenced by several uncertain factors such as volatile weather-dependent energy production, consumer behaviour, external market conditions, and waste availability. The resulting decisions can therefore be significantly impacted by uncertainty. This subsection outlines which aspects of uncertainty have been addressed in the reviewed literature by categorising the types of uncertainties considered. Uncertainty representation also includes how a specific uncertainty factor is implemented in a model e.g. through scenarios or uncertainty sets. For the sake of brevity, this section focuses solely on the different uncertainty factors, while modelling choices of uncertainty are presented alongside the methodology in Section \ref{Sec:Modelling}. An overview of this is given in Table \ref{tab:uncertainty}, which shows that 23 out of the 52 reviewed studies incorporate uncertainty aspects in their energy management problems.

\begin{table*}[h!]\centering
\footnotesize
\rowcolors{2}{gray!20}{white} 
\begin{tabular}{c c c c c l}
\multicolumn{4}{c}{\textcolor{foo}{\textbf{Uncertainty elements}}} & \phantom{abc} & \multicolumn{1}{c}{\textcolor{foo}{\textbf{Papers}}} \\
\cmidrule{1-4} \cmidrule{6-6}
Generation & Demand & Market Price & Waste Quality && \multicolumn{1}{c}{23}\\ \cmidrule{1-6}
X &&&&& \cite{WANG2023126288}, \cite{HWANGBO2022122006}, \cite{ZHAO2023101340}, \cite{XU2024123525}, \cite{CAI2025136178}, \cite{MOHTAVIPOUR2024132880}, \cite{chen_iet_nash}  \\
X & X &&&& \cite{LIU2024118204}, \cite{ALABI2022118997}, \cite{Zhang_DTU}, \cite{REN2023119499}, \cite{LI2022105261}, \cite{WEN2025134696}, \cite{ZHENG2025101098} \\
X & X & X &&& \cite{WANG2023121141}*, \cite{TOSTADOVELIZ2024123389}, \cite{ZHANG2023127643}, \cite{KHALIGH2023117354}, \cite{SALYANI2025106245} \\
X && X &&& \cite{SOHRABITABAR2022134898}, \cite{WANG20234631}\\
&& X &&& \cite{GUO202133039}\\
&& X & X && \cite{MOSTAFAVISANI2023116965} \\
\cmidrule{1-6}
\end{tabular}
\caption{Applied uncertainty representations in the reviewed literature, *: uncertainty wrt.\ EV arrival time is included}
\label{tab:uncertainty}
\end{table*}

Uncertainty regarding \textit{energy generation} is considered in 21 publications and addresses wind power \cite{WANG2023126288, Zhang_DTU,WANG2023121141, MOHTAVIPOUR2024132880, ZHENG2025101098}, solar power \cite{ZHAO2023101340, ALABI2022118997,LIU2024118204,TOSTADOVELIZ2024123389, LI2022105261} or both \cite{HWANGBO2022122006,XU2024123525, ZHANG2023127643,KHALIGH2023117354,SOHRABITABAR2022134898,WANG20234631, CAI2025136178, chen_iet_nash, WEN2025134696, SALYANI2025106245, REN2023119499}. \citet{Zhang_DTU} also include solar thermal production. 


In addition to energy generation uncertainty, \textit{energy demand uncertainty} is included by several studies. In many studies the uncertainty relates to the   power demand \cite{LIU2024118204, ALABI2022118997, Zhang_DTU, REN2023119499, LI2022105261, ZHANG2023127643, TOSTADOVELIZ2024123389, WANG2023121141, SALYANI2025106245, WEN2025134696, ZHENG2025101098}. 
Furthermore, uncertain thermal demands for heating \cite{SALYANI2025106245, WEN2025134696} both heating and cooling \cite{LIU2024118204, ALABI2022118997} or both heating and steam \cite{Zhang_DTU} are modelled. Finally,  hydrogen demand \cite{REN2023119499, TOSTADOVELIZ2024123389, SALYANI2025106245} or demand for EV charging \cite{WANG2023121141, SALYANI2025106245} are considered as uncertain.


Some studies further extend uncertainty representation by incorporating \textit{market price variability}. 
Most studies include uncertainty on electricity prices \cite{ZHANG2023127643,TOSTADOVELIZ2024123389,KHALIGH2023117354,WANG2023121141,SOHRABITABAR2022134898,WANG20234631,GUO202133039}. \citet{SALYANI2025106245} consider the balancing market. Other studies include also price uncertainty for other energy carriers such as hydrogen \cite{KHALIGH2023117354, MOSTAFAVISANI2023116965}, natural gas \cite{KHALIGH2023117354} or heat \cite{MOSTAFAVISANI2023116965}. The system in \citet{SOHRABITABAR2022134898} also sells cryptocurrency, so they also include cryptocurrency price uncertainty. 




The last uncertain factor is \textit{waste quality}, which is of particular importance to H4C due to relation with IS. 
The only study of the reviewed literature, that also account for waste quality uncertainty are \citet{MOSTAFAVISANI2023116965}, who consider multiple municipal solid waste compositions in their investment problem.

Overall, 23 of the 52 reviewed studies include uncertainty in their energy management models. These factors primarily focus on traditional uncertainty aspects such as generation, demand, and market prices. \unsure{Moreover, most of the considered uncertainties are related to the power sector. This suggests that the uncertainty frameworks applied in energy management models for H4C-related concepts are largely derived from models originating in the power sector.}

However, given the extensive waste sharing through IS in H4C, uncertainty in waste quality and availability may also significantly influence the optimisation decisions. Despite this, only one of the reviewed studies addresses uncertainty related to waste, highlighting a research gap in understanding how waste uncertainty affects energy management decisions within circular systems such as H4C. 

Ultimately, the choice of which uncertainties to include in an H4C energy management model should reflect the main uncertainties affecting the actual hub to be modelled. Including additional uncertainties that have minimal impact on decision-making may unnecessarily increase the computational complexity of the optimisation without significantly adding value to the optimised operation.

%% file: 3_5_operational_response.tex
\subsection{Operational Flexibility}\label{sec:flex}

To ensure a continuos energy supply, the production and consumption in an energy system have to be in balance at all times. Ensuring this balance is the primary purpose of energy management models, although the uncertainty factors discussed in Section \ref{sec:uncertainty} further complicate this task. To mitigate this, H4C energy management can leverage operational flexibility to schedule processes in response to various uncertainty realizations. Flexibility plays an important role in enabling IS operations, acting as an operational link between waste providers and receivers, and thereby driving the efficient utilization of waste resources. Operational flexibility can be broadly categorized into production, demand and storage flexibility. Historically, due to the perceived inflexibility of energy demand, energy production has been treated as the primary flexible component—scheduled to match relatively static demand profiles. However, the increasing integration of inflexible renewable energy sources has reversed this dynamic, creating a growing need for demand and storage flexibility \cite{White_sector}. Demand flexibility is achieved by shifting consumption in response to incentive signals such as dynamic pricing or emissions targets. Traditionally, this is referred to as demand response in the energy management literature. In the context of H4C, however, entities can both be energy consumers and waste suppliers through IS relationships. To reflect this duality, we adopt the broader term operational flexibility, which encompasses the adaptability of hub operations to support both IS efficiency and renewable energy integration. While both production and demand flexibility aim to synchronize production and consumption in time, storage introduces an additional degree of freedom by decoupling the two. 

The rest of this section presents our findings on applied operational flexibility in terms of demand response and storage solutions, covered in Sections \ref{sec:DR} and \ref{sec:storage}, respectively. The review also identifies key research gaps within each flexibility type.

\subsubsection{Operational Response} \label{sec:DR} 

To cope with fluctuations in energy generation, market prices, or waste availability, flexible consumers can adjust their consumption, either implicitly by responding to incentive signals or explicitly by offering their flexibility. Thereby, they can reduce the uncertainty effect on the energy management. This section focus on how this operational response has been applied throughout the reviewed literature. Table \ref{tab:DR} provides an overview of the studies including operational response and whether implicit or explicit operational response mechanisms has been applied.
\edit{All operational responses identified in the reviewed studies were implemented as demand response. Therefore, the term demand response will be used when referring to these implementations.}

\begin{table*}[h!]\centering
\footnotesize
\rowcolors{2}{white}{gray!20}
\begin{tabular}{l c l} 
{\textcolor{foo}{\textbf{Operational response type}}} & \phantom{abc} & {\textcolor{foo}{\textbf{Papers}}} \\
\cmidrule{1-3}
Implicit && \cite{GHASEMIMARZBALI2023120351}, \cite{SOHRABITABAR2022134898}, \cite{Dwijendra}, \cite{GUO202133039}, \cite{LI2024119946}, \cite{HOU2024130617}, \cite{CHEN2024120734}, \cite{MUGHEES2023121150}, \cite{MOHTAVIPOUR2024132880}, \cite{ZUO2025135413}, \cite{ZHENG2025101098} \\
Explicit && \cite{WANG2023121141}, \cite{Chen_Wei}, \cite{QIAN2024110504}, \cite{SALYANI2025106245}, \cite{PU_10637267}, \cite{chen_iet_nash}, \cite{WEN2025134696} \\
\cmidrule{1-3}
\end{tabular}
\caption{Applied operational response types in the reviewed literature}
\label{tab:DR}
\end{table*}

In \textit{implicit demand response}, energy consumption reacts to price signals to improve the objective function. This approach is commonly known as price-based demand response \cite{ALBADI20081989, FREIREBARCELO2022107953}. In our SLR, 11 of the studies implement implicit demand response.  \citet{Dwijendra} implement it by optimizing the entire system load with regards to power supply price signals.  However, in most studies, not all consumption is considered flexible. As an example, \citet{GUO202133039}, \citet{SOHRABITABAR2022134898}, \citet{MOHTAVIPOUR2024132880} and \citet{ZUO2025135413} consider only a subset of the energy consumption to be flexible, while the remaining portion follows a fixed demand structure. To prevent over- or under-production, \citet{GUO202133039} and \citet{SOHRABITABAR2022134898} both add a constraint ensuring that the sum of the demand response over a given time period equals zero. \citet{GHASEMIMARZBALI2023120351} further distinguish between flexible and non-flexible loads, restricting flexible demand shifts to a limited time window of nearby hours. \citet{CHEN2024120734} expand this differentiation by dividing cooling consumption into reducible and fixed loads, while treating electricity consumption similarly to \citet{GHASEMIMARZBALI2023120351}.

Both \citet{HOU2024130617} and \citet{LI2024119946} adopt a different demand response approach based on a Stackelberg game. In this model, energy producers (acting as leaders) send price signals, which energy users can respond to in order to improve their objectives. Thus, demand response is determined as part of the regular market-clearing process. Similarly, Mughees et al. \cite{MUGHEES2023121150} establish a market in which energy producers and consumers can exchange energy to improve the model's objective, but without relying on a Stackelberg game. 

\edit{Similar to \citet{Dwijendra}, \citet{ZHENG2025101098} also consider the entire load as flexible. However, the implicit price signal reflects a combination of each load’s individual benefit from providing demand response. The signal is based on each load's price elasticity and the negatively weighted average energy benefit of demand response at the system level. Furthermore, a penalty is given for large deviations from originally planned consumption.}

In \textit{explicit demand response}, energy users are paid specific remuneration or incentive for providing demand response \cite{ALBADI20081989,FREIREBARCELO2022107953}. \citet{WANG2023121141}, \citet{Chen_Wei}, \citet{SALYANI2025106245}, \citet{PU_10637267} and \citet{chen_iet_nash} all implement this using a fixed remuneration structure. \citet{WANG2023121141} apply explicit demand response similarly to \citet{GUO202133039} by treating a portion of demand as flexible while ensuring total energy consumption remains unchanged. \citet{Chen_Wei} adopt a demand response structure comparable to \citet{GHASEMIMARZBALI2023120351} and \citet{CHEN2024120734}. They categorize loads as curtailable, levelled, or transferable. \edit{\citet{SALYANI2025106245} differentiates between transferable electric loads and curtailable thermal loads, while \citet{PU_10637267} consider classical energy loads as curtailable and reducible and job scheduling loads at a maritime port as adjustable.} \citet{QIAN2024110504} classify loads as transferable, reducible, or interruptible with a staggered incentive structure, providing remuneration at varying price levels to the demand response participants. \edit{In addition to fixed price, \citet{chen_iet_nash} also base their demand response decisions of curtailable and transferable load on the consumer utility of changing their consumption. In contrast to the fixed price remuneration, \citet{WEN2025134696} implement a dynamic price compensation for transferable and interruptible loads. }

In conclusion, 18 of the 52 reviewed studies incorporate demand response mechanisms. As the share of renewable, weather-dependent energy continues to grow \cite{White_sector} and waste availability is influenced by upstream production patterns, flexible consumer behaviour will be important for effectively utilising limited energy and waste resources. Without such flexibility, industries may face the need for significant investments in energy and resource storage units, or risk energy and waste resource oversupply or shortage. Thus, demand response mechanisms are essential for achieving high resource efficiency in the context of H4C.

\subsubsection{Storages} \label{sec:storage}
As mentioned in Section \ref{sec:flex}, storages have the ability to decouple the time of production from the time of consumption of energy and resources. Considering the sector coupled energy carriers, a storage in one carrier can also help decoupling production and consumption for another downstream carrier. Table \ref{tab:storage} presents the energy and resource sectors in which the reviewed studies considered energy storages.

 \begin{table*}\centering
\footnotesize
\rowcolors{2}{gray!20}{white} 
\begin{tabularx}{\textwidth}{c c c c c c c c c X}
\rowcolor{white}
\multicolumn{8}{c}{\textcolor{foo}{\textbf{Storage Sectors}}} & \phantom{abc} & \multicolumn{1}{c}{\textcolor{foo}{\textbf{Papers}}} \\
\cmidrule{1-8} \cmidrule{10-10}
Electricity & Heat & Cooling & Gas & Hydrogen & Waste & Carbon & Water && \multicolumn{1}{c}{42}\\ \cmidrule{1-10}
X &&&&&&&&& \cite{GHASEMIMARZBALI2023120351}, \cite{HWANGBO2022122006}, \cite{LIN2021117741}, \cite{TOSTADOVELIZ2024123389}, \cite{fan10878712} \\
X & X &&&&&&&& \cite{WANG2023121141}, \cite{en16248041}, \cite{MUGHEES2023121150}, \cite{DAVOUDI2022107889}, \cite{ALABI2022118997}, \cite{su16125016}, \cite{LI2024129784}, \cite{REN2025135422}, \cite{ZHENG2025101098} \\
X & X & X &&&&&&& \cite{LI2024119946} \\
X & X & X & X &&&&&& \cite{SOHRABITABAR2022134898} \\
X & X &   & X & X &&&&& \cite{GUO202133039}, \cite{SALYANI2025106245}, \cite{chen_iet_nash} \\
X & X &   & X &   & X &&&& \cite{ZHAO2023101340} \\
X & X &   & X &&&&&& \cite{QIAN2024110504}, \cite{Chen_Wei}, \cite{LIU2024133528} \\
X & X &   &   & X &&&&& \cite{HOU2024130617}, \cite{MOSTAFAVISANI2023116965}, \cite{PANG2023125201}, \cite{CAI2025136178}, \cite{WEN2025134696} \\
X &   & X &&&&&&& \cite{WANG2024129868}, \cite{Zhang_DTU} \\
X &   &   &   & X &   & X&&& \cite{KHALIGH2023117354} \\
X &   &   &   & X &&&&& \cite{GOH2023117484}, \cite{WANG20234631}, \cite{WEI2022121732}, \cite{PU_10637267} \\
  & X & X &   & X &&&&& \cite{ZHANG2023127643}, \cite{DONG20251122}, \cite{ZUO2025135413} \\
  & X &   & X & X &&&&& \cite{REN2023119499} \\
  & X &&&&&&&& \cite{CHINESE2022124785} \\
  & X &   &   &   &   &  & X && \cite{MALDET2024100106} \\
  &   &   & X &   &   & X &&& \cite{LI2022105261} \\
\cmidrule{1-10}
\end{tabularx}
\caption{Considered storage mediums in the reviewed literature}
\label{tab:storage}
\end{table*}

Table \ref{tab:storage} shows that 42 out of 52 studies include storages in their mathematical problems. 35 of these include electricity storages in the form of batteries \cite{GHASEMIMARZBALI2023120351, HWANGBO2022122006, LIN2021117741, TOSTADOVELIZ2024123389, WANG2023121141, en16248041, MUGHEES2023121150, DAVOUDI2022107889, ALABI2022118997, su16125016, LI2024129784, LI2024119946, SOHRABITABAR2022134898, GUO202133039, ZHAO2023101340, QIAN2024110504, Chen_Wei, HOU2024130617, MOSTAFAVISANI2023116965, PANG2023125201, WANG2024129868, Zhang_DTU, KHALIGH2023117354, GOH2023117484, WANG20234631, WEI2022121732, fan10878712, CAI2025136178, PU_10637267, LIU2024133528, REN2025135422, WEN2025134696, ZHENG2025101098, SALYANI2025106245, chen_iet_nash}. This might be to mitigate the fluctuating renewable supply and due to the required balancing of energy supply and consumption at all times. The second largest storage sector is heat, which 29 studies include \cite{WANG2023121141, en16248041, MUGHEES2023121150, DAVOUDI2022107889, ALABI2022118997, su16125016, LI2024129784, LI2024119946, SOHRABITABAR2022134898, GUO202133039, ZHAO2023101340, QIAN2024110504, Chen_Wei, HOU2024130617, MOSTAFAVISANI2023116965, PANG2023125201, ZHANG2023127643, REN2023119499, CHINESE2022124785, MALDET2024100106, CAI2025136178, LIU2024133528, DONG20251122, REN2025135422, WEN2025134696, ZUO2025135413, ZHENG2025101098, SALYANI2025106245, chen_iet_nash}. Seven studies include cooling storage\cite{LI2024119946, SOHRABITABAR2022134898, WANG2024129868, Zhang_DTU, ZHANG2023127643, DONG20251122, ZUO2025135413}, while ten studies optimise systems with gas storages \cite{SOHRABITABAR2022134898, GUO202133039, ZHAO2023101340, QIAN2024110504, Chen_Wei, REN2023119499, LI2022105261, LIU2024133528, SALYANI2025106245, chen_iet_nash}. 17 studies consider hydrogen storages \cite{GUO202133039, HOU2024130617, MOSTAFAVISANI2023116965,PANG2023125201, KHALIGH2023117354, GOH2023117484, WANG20234631, WEI2022121732, ZHANG2023127643, REN2023119499, CAI2025136178, PU_10637267, DONG20251122, WEN2025134696, ZUO2025135413, SALYANI2025106245, chen_iet_nash}. As the only study, \citet{ZHAO2023101340} consider an urban and rural waste storage. \citet{KHALIGH2023117354, LI2022105261} store and utilise carbon in their studies, while \citet{MALDET2024100106} include a water storage. 

Overall, the review showed that most of the optimised systems include multiple types of storage, with 37 studies featuring storage across multiple sectors. This cross-sector integration enhances system flexibility compared to systems that include storage in only a single sector. Furthermore, only four studies incorporated storage within non-energy sectors, indicating that most energy management research focuses primarily on leveraging energy system flexibility. However, non-energy-based IS also require a balance between waste providers and receivers, where storage can offer the necessary flexibility to support this exchange. Therefore, investigating how non-energy storage influences energy management in H4C presents a relevant gap for future research.

%% file: 3_6_energy_market.tex
\subsection{Market Participation} \label{sec:market}
External markets can also act as an uncertainty mitigator by enhancing the supply resilience of the hub. In addition, external markets can also serve as an additional revenue stream. In contrast, internal markets function as platforms for resource and energy exchange within the hub itself, supporting local resource consumption. Alternatively, a hub may operate as an isolated system consuming only locally produced energy and resources. Therefore, another important modelling consideration is whether energy markets should be included, and if so, how they should be integrated, as this will significantly influence the resulting optimization strategies. This section outlines how markets have been applied throughout the reviewed literature. Table \ref{tab:Market} provides an overview of this, including whether a given study implements an internal market.

\begin{table*}\centering
\footnotesize
\rowcolors{2}{gray!20}{white} 
\begin{tabularx}{\textwidth}{c c c c c c c c X}
\rowcolor{white}
\multicolumn{7}{c}{\textcolor{foo}{\textbf{Energy Market}}} & \phantom{abc} & \multicolumn{1}{c}{\textcolor{foo}{\textbf{Papers}}} \\
\cmidrule{1-7} \cmidrule{9-9}
Electricity & Gas & Heat & Hydrogen & PtX-fuel & Cooling & Carbon && \multicolumn{1}{c}{47}\\ \cmidrule{1-9}
X &&&&&&&&  \cite{LIU2024118204}, \cite{MALDET2024100106}, \cite{SOHRABITABAR2022134898}, \cite{Chen_Wei}, \cite{en16248041},  \cite{Zhang_DTU}, \cite{CAI2025136178}, \cite{PU_10637267}$^M$, \cite{ZHENG2025101098}, \cite{fan10878712} \\
X & X &&&&&&& \cite{WANG2023126288}$^M$, \cite{MENEGHETTI2012263}, \cite{WEI2022121732}, \cite{DAVOUDI2022107889}, \cite{LI2024119946}$^M$, \cite{LIN2021117741}, \cite{CHINESE2022124785}, \cite{QIAN2024110504}, \cite{LI2024129784}, \cite{su16125016}, \cite{MOHTAVIPOUR2024132880}, \cite{REN2025135422}, \cite{chen_iet_nash}$^M$, \cite{pedro10161834}, \cite{ZUO2025135413}$^M$  \\
X & X & X&&&&&& \cite{WANG2023121141}, \cite{Dwijendra}, \cite{GUO202133039}, \cite{WANG2024129868}$^M$, \cite{ZHAO2023101340}, \cite{LIU2024133528}\\ 
X & X && X&&&&&  \cite{REN2023119499}$^M$, \cite{SALYANI2025106245} \\
X &   & X & X &&&&& \cite{MOSTAFAVISANI2023116965} \\
X & & & X &&&&& \cite{TOSTADOVELIZ2024123389}, \cite{ZHANG2023127643}, \cite{XU2024123525}, \cite{WANG20234631}, \cite{GOH2023117484}, \cite{DONG20251122} \\
X & X && X& X&&&& \cite{KHALIGH2023117354} \\
X && X &&&&&& \cite{MUGHEES2023121150}$^M$\\
X &&&&& X&&& \cite{CHEN2024120734}$^M$ \\
X & X &&&&& X && \cite{HOU2024130617}$^M$, \cite{WEN2025134696}$^M$ \\
& X&&&&&&& \cite{ALABI2022118997}, \cite{LI2022105261} \\
\cmidrule{1-9}
\end{tabularx}
\caption{Interaction with markets in hte reviewed literature, \textbf{M}: internal hub market}
\label{tab:Market}
\end{table*}

Table \ref{tab:Market} shows that 47 out of the 52 reviewed studies include a form of market in their energy management model. Furthermore, Table \ref{tab:Market} also shows that the electricity market is the main considered market, 45 studies including it in their models. Most studies solely considering the electricity market \cite{LIU2024118204,MALDET2024100106,SOHRABITABAR2022134898,Chen_Wei,en16248041,Zhang_DTU,CAI2025136178,ZHENG2025101098,fan10878712} or gas market \cite{ALABI2022118997,LI2022105261}, and all assume a price-taker perspective. Similarly, most studies considering combinations of several energy markets consider a price-taker approach. The studied combinations are: electricity and gas \cite{MENEGHETTI2012263,WEI2022121732,DAVOUDI2022107889,LIN2021117741,CHINESE2022124785,QIAN2024110504,LI2024129784,su16125016,MOHTAVIPOUR2024132880,REN2025135422,pedro10161834}; electricity, gas and heat \cite{WANG2023121141,Dwijendra,GUO202133039,ZHAO2023101340,LIU2024133528}; electricity (both day-ahead and balancing), gas and hydrogen \cite{SALYANI2025106245}; electricity heat and hydrogen \cite{MOSTAFAVISANI2023116965}; electricity and hydrogen \cite{TOSTADOVELIZ2024123389,ZHANG2023127643,XU2024123525,WANG20234631,GOH2023117484,DONG20251122}. Khaligh et al. \cite{KHALIGH2023117354} further extends this by also incorporating the markets for gas and hydrogen byproducts markets in the form of urea and methanol (PtX-fuel) in their study. \edit{In contrast to the previously mentioned studies only considering energy markets, \citet{WEN2025134696} also incorporate a carbon market along with an electricity market and a gas market.}

Apart from external markets, \textit{internal markets}, where market prices are outcomes of the optimisation, are also used throughout the literature. \citet{PU_10637267} and \citet{chen_iet_nash} establish internal markets for electricity, while \citet{WANG2023126288}, \citet{LI2024119946}, \citet{ZUO2025135413} and \citet{HOU2024130617} model combined internal gas and electricity markets. In their study, \citet{WANG2023126288} investigates how a small portfolio of strategic energy suppliers can improve their profit by bidding strategically in the internal market and how this would alter the market outcome. \citet{PU_10637267}, \citet{chen_iet_nash}, \citet{LI2024119946}, \citet{HOU2024130617} and \citet{ZUO2025135413} all establish a market clearing problem, where an energy operator formulates prices that the energy consumers react to by changing their consumption. Furthermore, \citet{HOU2024130617} also include a carbon emission quota and green certificate scheme in their model. 
\citet{MUGHEES2023121150} model their own internal market in which electricity and heat can be traded to improve the internal energy efficiency. Furthermore, if needed, external electricity can be purchased from the traditional electricity market.  \citet{CHEN2024120734} introduce an energy retailer that acts as an intermediary between suppliers and consumers, managing the scheduling of electricity and cooling. In addition to the external market, \citet{LI2022105261} further incorporates an energy retailer that takes care of the energy sharing internally in the hub. In \citet{WANG2024129868} and \citet{CAI2025136178}, multiple energy hubs trade collaboratively before trading with the external markets to be locally cost-efficient. \citet{REN2023119499} include external electricity, gas and hydrogen markets, which they model sequentially similar to \citet{LI2024119946} and \citet{HOU2024130617}.



Overall, market participation in the reviewed literature has been limited to either internal market implementations or price-taker assumptions. Both approaches rely on the presence of a central entity, where hub participants either allow the entity to manage an internal market or to participate in external markets on behalf of the entire hub. This effectively treats the hub as a portfolio. However, since concepts such as H4C involve various companies collaborating through both bottom-up and top-down initiated industrial symbiosis (IS), the presence of a focal company or hub leader capable of assuming this role cannot necessarily be assumed. In some cases, the collaboration may not be sufficiently strong to support such a centralised structure. Consequently, this portfolio optimisation approach may not be applicable to all types of H4C. This will be further discussed in Section \ref{Sec:Modelling}. 

Furthermore, the external market participation in the reviewed studies largely relies on a price-taker assumption. None of the studies actively model market bidding, meaning they exclude the sequential dynamics \edit{resulting from determining and submitting bids to markets and market clearance} in their optimisations. Given the complex interconnectivity between multiple energy and resource carriers in concepts such as H4C, this presents a gap for future research. Additionally, in H4C, waste markets could facilitate waste sharing among hub processes. However, as only two studies consider how to re-utilise waste in the form of carbon in their models, this shows another gap within the literature.

%% file: 3_7_Centralised.tex
\subsection{Approaches for Collaboration and Information Sharing} \label{Sec:Modelling}

Due to the IS dynamics and sector-coupled energy flows in H4C, resource and energy sharing is becoming increasingly decentralised \cite{White_sector}. This results in multiple actors being operationally reliant on the supply and, thereby, the operation of other actors in the hub. To coordinate this complex energy and waste sharing, information exchange and collaboration are needed. However, this has historically proven to be a challenge in facilitating industrial symbiosis \cite{PATALA2020121093,SONG2018414}. IS relationships can both be initiated from top down facilitation by a central mediator and from a bottom up initiative between two individual hub processes. In a hub primarily based on bottom up initiated IS relationships, an inherent natural central entity might not be present. In these kind of decentralised systems, the hub actors might not be willing to share information with a central entity, and therefore the information should only be shared among the actors reliant on each other. Secondly, it is not necessarily given that they are willing to cooperate with the whole hub. \citet{JEONG2023106190} argues that an appropriate optimisation model for such a decentralised supply should also resemble this decentralised structure. 
Conversely, in hubs where IS is primarily initiated top-down, a central entity facilitates collaboration, making actors more inclined to share information and collaborate with all hub participants, even those with whom they have no direct IS relationship. If this is the case, a centralised model might be appropriate to represent the actual hub dynamics.

This section categorises the reviewed articles into centralised and decentralised modelling approaches and outline the various applied modelling methodologies. Furthermore, we discuss the strengths and limitations of these approaches while identifying gaps in the reviewed literature.

\subsubsection{Centralised Approaches} \label{sec:central}

Among the 52 studies included in the SLR, 47 incorporated a central operator in their energy management models. Table \ref{tab:centralised} categorises the different studies according to their applied modelling methodologies.

Energy hubs are typically modelled by distinguishing between different energy carriers and resource types within a hub, then employing conversion factors to represent the various conversion and production processes occurring throughout the hub \cite{MOHAMMADI20171512}.  Common modelling approaches for this are using a deterministic \textit{linear program }(LP) \cite{pedro10161834, fan10878712} or a deterministic \textit{ mixed-integer linear program }(MILP) \cite{MALDET2024100106,MENEGHETTI2012263,WEI2022121732,Chen_Wei,en16248041,CHINESE2022124785, Kantor_Frontier, LIN2021117741, LIU2024133528}. 

When uncertainty is introduced into the energy management process, various modelling approaches are applied, differing primarily in how they represent and handle that uncertainty. Stochastic and chance constraint programming treat uncertainty through multiple scenarios with assigned outcome probabilities, while robust optimisation optimises for the worst-case realisation. 

In \textit{stochastic programming}, the decisions are divided into stages, where decisions made in a previous stage affect the decisions in later stages. The optimisation aims to find the best objective while ensuring feasibility across all scenarios. The most common objective is the expected value. \citet{LIU2024118204} and \citet{KHALIGH2023117354} apply the expected value using an MILP, while \citet{XU2024123525} apply an MINLP. Furthermore, \citet{MOHTAVIPOUR2024132880} apply the expected value approach to an originally non-convex problem, which they  relax to a semidefinite program that is easier to solve. However, the expected value approach might result in favourable outcomes at the expense of poorer-performing scenarios. To mitigate this, a term prioritising the optimization of poor-performing scenarios can be added to the optimisation, such as including the conditional value at risk (CVaR). \citet{SOHRABITABAR2022134898}, \citet{WANG20234631} and \citet{HWANGBO2022122006} all apply CVaR using an MILP. However, in the study by \citet{HWANGBO2022122006}, the stochastic solution is only a part of a broader optimisation that first optimise a deterministic program to determine non-renewable capacity installations before using the stochastic solution to determine the renewable capacity mix. 
\textit{Chance constraint programming}, which \citet{ZHANG2023127643} apply, similarly represents uncertainty through scenarios but permits certain constraints to be violated for a limited number of cases, promoting more opportunistic decision-making.

When more risk-adverse behaviour is preferred, \textit{robust optimisation approaches} can be applied. \citet{GUO202133039} apply classical robust optimisation, optimising for the worst-case realisation within an uncertainty set. On the other hand, \citet{MOSTAFAVISANI2023116965} consider a problem with adaptive robust optimization with two stages like the stochastic program. Zhao et al. \cite{ZHAO2023101340} propose a distributionally robust optimisation model that consider an ambiguity set of possible distribution set for the uncertain parameters. They compute the optimal decisions based on the worst uncertainty distribution within the set.



Unlike the previously mentioned robust optimisation programs, \citet{ALABI2022118997}, \citet{TOSTADOVELIZ2024123389}, \citet{SALYANI2025106245}, \citet{LI2022105261}, \citet{WEN2025134696} apply a combination of stochastic programming and robust optimisation. \citet{ALABI2022118997} combine scenarios for the energy generation with a robust representation of load uncertainty using duality theory in their MILP. Similarly, \citet{TOSTADOVELIZ2024123389} consider stochastic electricity prices but apply a robust formulation for the PV production and electricity demand.  They solve it in multiple stages handling investment, binary, and operational decisions separately using Karush–Kuhn–Tucker conditions. \citet{SALYANI2025106245} consider a robust program to handle real time market prices and stochastic renewable production and demands. The prices are handled using a budget of uncertainty, while they use the CVaR method in combination with probability density functions to handle the stochasticity of generation and demand. \citet{LI2022105261} model a stochastic program, where they apply a P-robust measure ensuring that an objective value for each scenario cannot be worse than a given factor of the deterministic optimal result for the given scenario. \citet{WEN2025134696} also intend to obtain a less conservative outcome by combining robust optimisation through uncertainty set with chance constraints. 

The structure in \citet{Dwijendra}, \citet{QIAN2024110504} and \citet{Zhang_DTU} is also divided but in \textit{multiple optimisation layers} that are solved iteratively. Both \citet{Dwijendra} and \citet{QIAN2024110504} utilise the multi-layer structure to optimise demand response, which they both solve using different heuristics. Meanwhile,  \citet{Zhang_DTU} solve a stochastic program as upper layer for a full hour, while the second operational layer determine the actual operation based on a more process-detailed program running at a shorter time horizon. 
Similarly, \citet{PANG2023125201} and \citet{GHASEMIMARZBALI2023120351} use an MINLP to capture the complex process dynamics. \citet{PANG2023125201} use non-linear terms to represent solid oxide fuel and electrolyser cell processes more detailed while \citet{GHASEMIMARZBALI2023120351} apply a quadratic program to include the dynamic inlet valve effects in the the relationship between fuel costs and thermal power output. Such non-linear terms may reduce computational tractability compared to MILP approaches.

\begin{table*}\centering
\footnotesize
\rowcolors{2}{white}{gray!20} 
\begin{tabular}[h!]{l c l}
{\textcolor{foo}{\textbf{Centralised modelling}}} & \phantom{abc} & {\textcolor{foo}{\textbf{Papers}}} \\
\cmidrule{1-3}
Deterministic (MI)LP && \cite{pedro10161834}, \cite{fan10878712}, \cite{CHINESE2022124785},  \cite{LIN2021117741}, \cite{MENEGHETTI2012263},  \cite{Kantor_Frontier}, \cite{Chen_Wei}, \cite{en16248041}, \cite{WEI2022121732}, \cite{MALDET2024100106} \\
Stochastic programming && \cite{KHALIGH2023117354}, \cite{WANG20234631}, \cite{SOHRABITABAR2022134898},   \cite{XU2024123525}, \cite{HWANGBO2022122006}, \cite{LIU2024118204}, \cite{LIU2024133528}, \cite{MOHTAVIPOUR2024132880} \\
Chance constraint programming && \cite{ZHANG2023127643} \\
Robust optimisation &&  \cite{ZHAO2023101340}, \cite{LI2022105261}, \cite{MOSTAFAVISANI2023116965}, \cite{ALABI2022118997}, \cite{GUO202133039}, \cite{TOSTADOVELIZ2024123389}, \cite{SALYANI2025106245}, \cite{WEN2025134696} \\
Multi-layer MILP && \cite{Zhang_DTU}, \cite{Dwijendra}, \cite{QIAN2024110504} \\
(MI)NLP && \cite{GHASEMIMARZBALI2023120351}, \cite{PANG2023125201} \\
Stackelberg Game && \cite{LI2024119946}, \cite{WANG2023126288}, \cite{CHEN2024120734}, \cite{GALVANCARA2022117734}, \cite{HOU2024130617}, \cite{chen_iet_nash}, \cite{ZUO2025135413} \\
Heuristic && \cite{REN2023119499}, \cite{GOH2023117484}, \cite{WANG2023121141}, \cite{DAVOUDI2022107889}, \cite{LI2024129784}, \cite{REN2025135422}, \cite{ZHENG2025101098} \\
Reinforcement Learning && \cite{MUGHEES2023121150} \\
\cmidrule{1-3}
\end{tabular}
\caption{Centralised modelling approach overview}
\label{tab:centralised}
\end{table*}

An alternative modelling approach is to establish an internal hub energy sharing market. This is typically done by modelling the strategic \textit{Stackelberg game}, where a subset of the participants, classified as leaders, initially make a decision in an upper level program, while the remaining participants, classified as followers, react to the leaders' decisions in a lower level program. This negotiation wants to reach an optimisation point, where both the leaders and followers do not want to change their decision, which is the Nash equilibrium. \citet{HOU2024130617}, \citet{LI2024119946}, and \citet{GALVANCARA2022117734} group the energy producers or a hub operator as leaders, while they model classical energy loads as followers. In these games, the leaders formulate prices to facilitate a demand response from the followers. On the other hand, \citet{CHEN2024120734} consider an energy retailer formulating dynamic prices as leader, which the energy suppliers and consumers respond to as followers. \edit{\citet{chen_iet_nash} consider an operator of multiple microgrids as the leader, who formulates energy prices for the microgrid-followers. They solve it using an iterative search algorithm. \citet{ZUO2025135413} consider an operator for each connected energy carrier that has to agree on energy scheduling. To optimise this, they combine the operator optimisations using Karush–Kuhn–Tucker conditions and thereby mathematically creating a leader operator that schedules all the energy carriers.}  Lastly, in the Stackelberg game modelled by \citet{WANG2023126288}, they distinguish between strategic and non-strategic energy producers. Contrary to the other Stackelberg game studies, they model the followers' lower level using distributionally robust optimisation to react to the realisation of wind power uncertainty.

To address computational complexity, \citet{DAVOUDI2022107889}, \citet{GOH2023117484}, \citet{LI2024129784}, \citet{REN2023119499}, \citet{LI2024129784}, \citet{REN2025135422} and \citet{ZHENG2025101098} apply heuristics to approximate pareto fronts for multi-objective problems. \citet{DAVOUDI2022107889} use a multi-objective particle swarm optimisation, while \citet{GOH2023117484} apply multi-objective mayfly algorithm in combination with the technique for order of preference by similarity to ideal solutions method. \citet{REN2023119499} apply an interior point method and a multi-objective golden eagle optimiser sequentially to solve their problem. \citet{LI2024129784} formulate their own gradient based optimizer in combination with fitness values and the local escaping operator method to handle the multiple objectives, while \citet{REN2025135422} use a multi-objective gradient-based optimizer with shift-based density estimation. Lastly, \citet{ZHENG2025101098} use the mountain
gazelle optimizer heuristic to both optimise across multiple objectives and optimise the demand response.  In contrast, \citet{WANG2023121141} uses a local escaping operator method in combination with a gray wolf algorithm to solve their robust problem. 
Finally, Mughees et al. \cite{MUGHEES2023121150} implement a Markov decision process solved via reinforcement learning. While this method may require extensive training time, it provides fast solutions once trained.

In conclusion, 47 out of the 52 reviewed studies implement a centralised modelling perspective for EM problems, with mathematical programming models based on LP and MILP being the most common. Only ten studies consider different methodologies such as MINLP, reinforcement learning or heuristics. 

\subsubsection{Decentralised Approaches} \label{sec:decentral}

In the SLR, only five studies included a decentralised modelling structure for their energy management problem. To provide a more comprehensive overview of the decentralised modelling techniques that are being applied throughout the literature, six additional studies within similar topics are included in the review for this section. These six studies were not captured in the original query, since they did not include the keywords related to circular economy (\textit{TS3}). However, despite not being part of an additional systematic literature search, they provide overall insights into how decentralised modelling approaches can be applied for coupled energy systems. 

A decentralised model structure does not necessarily imply the absence of a central hub operator. Throughout the literature, multiple decentralised energy management structures are applied, where a central operator is present either to handle energy transmissions or to facilitate the convergence of the mathematical optimisation process. Furthermore, the term \textit{decentralised energy management} does not imply that the scheduling of the energy actors is exclusively on an individual basis. In the existing literature, both cooperative and non-cooperative energy scheduling mechanisms are observed. To address this, Table \ref{tab:decentralised} categorises eleven studies based on the presence of a hub operator in the model and whether energy scheduling is modelled as cooperative or non-cooperative.

\begin{table*}[!h]\centering
\footnotesize
\rowcolors{2}{gray!20}{white} 
\begin{tabular}{l c c c}
{\textcolor{foo}{\textbf{Paper}}} & {\textcolor{foo}{\textbf{Methodology}}} & {\textcolor{foo}{\textbf{Central Operator}}} & {\textcolor{foo}{\textbf{Cooperative game}}} \\
\cmidrule{1-2} \cmidrule{3-4}
 \citet{WANG2024129868}               & ADMM                     & X & X \\
 \citet{WANG2020106060}               & ADMM                     & X &   \\ 
 \citet{ZAREIGOLAMBAHRI2024109901}    & ADMM                     & X & X \\
 \citet{HUSSAIN2024110458}            & ADMM                     &   &   \\
 \citet{DONG20251122}                 & Lagrangian Decomposition & X & X \\
\citet{PU_10637267}                   & Stackelberg              & X &   \\
\citet{CAI2025136178}                 & Stackelberg              & X &   \\
\citet{ZHOU2018993}                   & Nash-based algorithm     & X &   \\
\citet{CHEN2023128641}                & RL                       & X & X \\
 \citet{ZHU2022118636}                & RL                       &   & X \\
\citet{su16125016}                    & Heuristic                & X  & X \\
\cmidrule{1-4}
\end{tabular}
\caption{Decentralised modelling approach overview}
\label{tab:decentralised}
\end{table*}

To solve decentralised scheduling models, six of the included articles apply iterative game-theoretic based solutions to simulate a bargaining process. They iterate the scheduling process until a Nash Equilibrium is established.
\citet{WANG2024129868}, \citet{WANG2020106060}, \citet{HUSSAIN2024110458} and \citet{ZAREIGOLAMBAHRI2024109901} establish this by the \textit{alternating direction method of multiplier} (ADMM) algorithm. ADMM introduces a Lagrangian multiplier-penalty term together with an augmented quadratic penalty that forces the different agents to schedule towards an equilibrium \cite{boyd2011admm}. \citet{WANG2024129868} and \citet{WANG2020106060} use the ADMM algorithm to solve scheduling problems of multiple energy hubs trading energy with each other. Here, the central hub operators of multiple hubs trade energy to lower their costs. \citet{WANG2024129868} solve this through cooperative trading in two stages. In the first stage, each hub operator optimises their own asset portfolio, while all the hub operators trade together in the second stage. In the study by \citet{WANG2020106060}, they use the ADMM algorithm to solve a single-stage non-cooperative trading process. 
In \citet{ZAREIGOLAMBAHRI2024109901}, a central electricity and heat operator minimises the amount of energy waste within a network. The operator achieves this by offering flexibility deals through peer-to-peer trading with energy agents solved using ADMM.  To solve a demand response program for an energy community, Hussain et al. \cite{HUSSAIN2024110458} use ADMM with agents trading peer-to-peer. In their model, agents share information only with their neighbouring agents, thus eliminating the need for a central operator and ensuring information privacy. 

\edit{Similar to ADMM, \citet{DONG20251122} apply Lagrangian decomposition, where each player optimises their own system. After optimisation, an operator gains access to the optimized energy flows and, to facilitate convergence, sends a Lagrangian penalty term back. \citet{PU_10637267} and \citet{CAI2025136178} also make use of Lagrangian penalties to support convergence. However, they structure their optimisation as the Stackelberg game previously explained in Section \ref{sec:central}. The leader acts as an operator, and both papers model the followers as non-cooperative players. In addition, \citet{CAI2025136178} apply distributionally robust optimization and incorporate chance constraints to reduce the conservativeness of the results.}

\citet{ZHOU2018993} optimise the decisions for an energy retailer. Their model coordinates energy trading among multiple prosumer agents, who both produce and consume energy. They coordinate based on a non-cooperative algorithm aiming to reach a Nash equillibrium with limited information sharing among the optimising agents.

\citet{CHEN2023128641} and \citet{ZHU2022118636} formulate their models as \textit{Markov decision processes}, which they solve using reinforcement learning. \citet{CHEN2023128641} consider multiple hubs trading together with a reinforcement learning agent for each hub. The agents are trained centrally such that they cooperatively make decisions based on the interests of all the hubs combined. \citet{ZHU2022118636} also train multiple agents centrally, but their model consist of a single hub with an agent for each energy supplying asset. 

Lastly, \citet{su16125016} employ an iterative metaheuristic known as the Zebra Optimization Algorithm, where in each iteration, the operator focuses on optimizing the production process identified as the most suboptimal. The remaining assets are then optimised individually. This heuristic proceeds iteratively until convergence is achieved.

In conclusion, only five out of the 52 reviewed studies apply a decentralised model structure to solve their energy management problems. However, decentralised model structures are applied within the literature of similar fields, where game-theoretic concepts, reinforcement learning agents, and heuristics are used to solve the models. The literature includes studies both with and without a central operator, as well as those assuming cooperative or non-cooperative behavior among players.

\subsubsection{Centralised vs. Decentralised}
From an optimization perspective, a model with partial information can, at best, match the performance of a model with full information when solved to optimality \cite{Feigenbaum2007363}.
However, given that a centralised approach must consider a larger set of agents, it may result in computational tractability issues. Conversely, a decentralised approach distributes this computational burden among several computing agents, potentially mitigating these challenges. From a mathematical and computational perspective, the modelling choice of the two approaches becomes a trade-off of accuracy and tractability. 

However, it cannot be assumed that a centralised modelling approach is applicable to all types of H4C.  \citet{PATALA2020121093} identify data sharing concerns and the absence of demonstrated value as primary barriers to facilitating IS through an intermediary. Given that centralised energy management for H4C requires more extensive data sharing and greater autonomy than IS facilitation, these barriers are likely even more pronounced in the energy management context.

\citet{PATALA2020121093}  further suggest that fostering community relations within a hub can encourage a collective resource-sharing mindset, which may reduce the impact of these barriers. This mindset could similarly enhance the acceptance of a centralised energy management system for H4C. Hubs with stronger collaboration may be more inclined to accept a centralised decision-making entity managing energy on behalf of all stakeholders within the hub. Conversely, hubs with limited collaboration may resist such an approach. \edit{Furthermore, actors that are geographically distant from one another may be more reluctant to accept the involvement of a central coordinating entity compared to actors located in closer proximity. In the literature, geographical proximity is typically regarded as an enabler of IS \cite{TLEUKEN2025125324}.}

Ultimately, the suitability of a centralised or decentralised energy management approach depends on the specific characteristics of the hub being modelled. Despite this, the dominance of centralised approaches in the literature suggests a potential modelling bias. The fact that 47 out of 52 reviewed papers adopt a centralised approach indicates that this preference may stem more from established electricity modelling conventions than from considerations of collaboration dynamics within individual hubs. This lack of decentralised modelling approaches for similar concepts to H4C indicate a gap in the literature.

%% file: 3_2_Framework.tex
\section{Modelling Framework for Energy Management in H4C} \label{sec:framework}
\begin{figure}[!h]
    \centering
    \includegraphics[width=\linewidth]{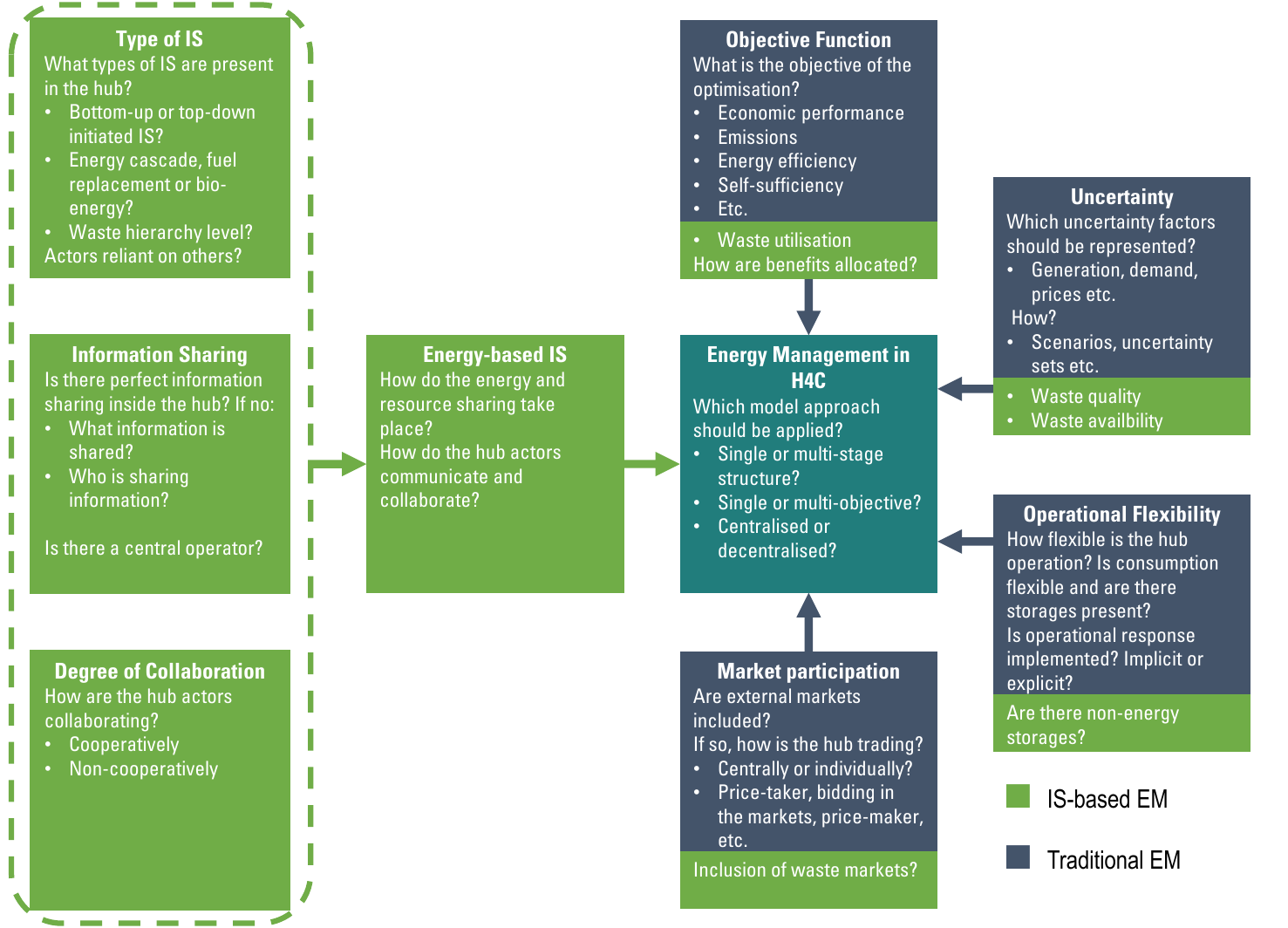}
    \caption{Modelling framework for energy management in hubs for circularity, EM: Energy Management}
    \label{fig:EM_in_H4C}
\end{figure}

In this section, we establish a framework to create an understanding of which aspects should be considered in an H4C energy management model. Based on the reviewed literature in Section \ref{sec:SLR}, we propose the framework illustrated in Figure \ref{fig:EM_in_H4C}.

We argue that energy management in H4C consists of both traditional energy management aspects, which are already being explored in the energy management literature (blue), and additional aspects arising from the industrial symbiosis (IS) taking place within H4C (green).

As in energy hubs, energy communities or EIPs, energy management methods in H4C need to include an objective function, handle uncertainties and leverage operational flexibility. Traditionally, energy management studies have prioritised economic performance in objective functions. However, with growing corporate interest in circularity \cite{wef2025circular}, more holistic objective functions incorporating emissions, energy efficiency, and other sustainability metrics are increasingly adopted. For H4C energy management, IS-related metrics such as waste utilisation may also be included in future objective functions. Furthermore, due to the involvement of multiple actors, the fair allocation of shared energy management benefits may need to be addressed.

To handle uncertainty, the first step is to identify the factors affecting hub operations. In traditional energy management, uncertainty primarily arises from renewable energy generation, demand fluctuations, and market price volatility. In waste-based IS within H4Cs, however, additional uncertainties may emerge from the availability and quality of waste resources. Once identified, these uncertainties need to be appropriately represented in the model, e.g. as scenarios in stochastic programs, as uncertainty sets in robust optimisation approaches or stochastic processes in Markov decision processes. \looseness=-1

To mitigate the impact of uncertainty on operation, the flexibility within a hub can be leveraged. Since renewable generation is weather-dependent and thus relatively inflexible \cite{White_sector}, flexibility can take the form of energy storage or the ability for consumers to shift their consumption or production. This allows operations to adapt to supply bottlenecks, e.g. low renewable production, and take advantage of more favourable conditions, e.g. high renewable production, when they arise. Moreover, in an H4C context, flexibility also plays a key role in facilitating IS operations, where especially storages between waste providers and receivers can provide the flexibility required for the IS relationship to operate. When modelling operational flexibility, key considerations include how adaptable hub operations are and how this flexibility can effectively be leveraged. \looseness=-1

In addition to operational flexibility, external markets may also serve to mitigate uncertainty by enhancing supply resilience through external procurement. Furthermore, markets can also be an additional revenue stream if surplus energy is sold to them. Alternatively, a hub may operate as an isolated system, relying solely on locally produced energy. Consequently, another important modelling consideration is whether and how energy markets should be included, as their integration significantly influences optimisation outcomes.

Beyond these traditional aspects, IS-specific factors must also be addressed when modelling energy management in H4C. Specifically, we argue that the type of IS within the hub, the degree of information sharing, and how the different hub processes collaborate are interconnected and should be explicitly incorporated into the modelling approach. As discussed in the introduction to Section \ref{Sec:Modelling}, both the structure of information exchange and the level of actor collaboration shape the suitability of different modelling approaches. These aspects also relate to the IS type: for example, a hub with primarily top-down IS relationships may more naturally involve a central coordinator, while a hub with bottom-up IS relationships may be more reluctant. \edit{Furthermore, if waste requires treatment before it can be used, such as at the recycling or incineration level in fuel replacement or bioenergy, or requires to be transported, additional actors may be involved. This can result in the operations of multiple actors becoming reliant on the waste output on upstream waste providers and treatment facilities. Additionally, from a circular economy perspective, energy flows on lower waste hierarchy levels need to be prioritised in energy management.}

In conclusion, the choice between centralised and decentralised modelling approaches depends on how IS relationships are structured, how information is shared, and how hub actors collaborate. These factors, together with considerations regarding uncertainty representation, market participation, and operational flexibility, define the overall decision-making structure of the energy management model for a given H4C.

%% file: 4_Conclusion.tex
\section{Conclusion} \label{sec:conclusion}

Based on a thematic systematic literature review of concepts similar to H4C, we proposed a modelling framework for energy management practitioners and researchers outlining key aspects that should be considered when developing H4C energy management models. We argue that energy management within H4C builds on established modelling concepts from interconnected energy systems but must also reflect the complex information and collaboration structures arising from the energy-based IS taking place within H4C. 

The review indicates a degree of modelling path dependency on traditional power system modelling. Like traditional power energy management, many studies rely on centralised modelling approaches and power-related uncertainties, such as generation and demand uncertainty. Our proposed framework aims to address this by guiding future modelling efforts to incorporate IS dynamics in order to more accurately reflect the nature of energy sharing within H4C. For this an extension of the traditional approaches is necessary, explicitly addressing the more or less extensive collaboration among several stakeholders as well as energy and waste flows apart from power, gas and heat.

As the H4C concept was only formally introduced in 2021, no H4C energy management literature has to the best of our knowledge been published at the time of article submission. Therefore, the conclusions drawn from similar concepts may not fully represent how H4C energy management models will evolve. It is also possible that additional aspects will need to be integrated into the framework over time. Future research should address gaps in decentralised energy management modelling and the treatment of waste-related uncertainties. Additionally, exploring the potential synergies between energy markets and IS-based energy sharing in H4C would provide valuable insights for future model development.

%% file: AppendixTable.tex
\appendix
\section{Appendix: Study Overview} \label{App:tab}
\newgeometry{left=1.5cm, right=1.5cm, top=2cm, bottom=2cm}
\begin{landscape}
\begin{table}[ht]
\centering
\resizebox{\linewidth}{!}{
\begin{tabular}{lccccccccc}
\textbf{Study} & \textbf{Concept} & \textbf{Waste Hierarchy} & \textbf{IS} & \textbf{Objective} & \textbf{Uncertainty} & \textbf{Operational response} & \textbf{Storages} & \textbf{Markets} & \textbf{Modelling Approach} \\
\midrule
\citet{ALABI2022118997}& Energy Hub & R2 & - & ECO, EMI, CUR \& DEG & GEN \& DEM & - & E \& HE & HE & Robust (CEN) \\
\citet{CAI2025136178}& Energy Hub & R2 \& R3 & EC & ECO \& EMI & GEN & - & E, HE \& HY & E & Stackelberg (DEC) \\
\citet{Chen_Wei}& Industrial Park & R2, R8 \& R9 & FR \& BE & ECO \& EMI & - & Explicit & E, HE \& G & E & Deterministic MILP (CEN) \\
\citet{chen_iet_nash}& Industrial Park & R2 \& R8 & - & ECO, EMI, CUR \& DRU & GEN & Explicit & E, HE, G \& HY & E \& G & Stackelberg (CEN) \\
\citet{CHEN2024120734}& Energy Hub & R2 & - & ECO \& EMI & - & Implicit & - & E \& CO$^M$ & Stackelberg (CEN) \\
\citet{CHINESE2022124785}& Energy Hub & R2 \& R3 & EC & ECO \& EMI & - & - & HE & E \& G & Deterministic MILP (CEN) \\
\citet{DAVOUDI2022107889}& Energy Hub & R2 & - & ECO \& ENE & - & - & E \& HE & E \& G & Heuristic (CEN) \\
\citet{DONG20251122}& Industrial Park & R2 & - & ECO \& EMI & - & - & HE, CO \& HY & E \& HY & Lagrangian Decomposition (DEC) \\
\citet{Dwijendra}& Energy Hub & R2 & - & ECO & - & Implicit & - & E, G \& HE & Multi-layer MILP (CEN) \\
\citet{fan10878712}& Energy Hub & R2 \& R8 & BE & ECO \& EMI & - & - & E & E & Deterministic LP (CEN) \\
\citet{GALVANCARA2022117734}& Industrial Park & R2 \& R3 & EC & ECO \& EMI & - & - & - & - & Stackelberg (CEN) \\
\citet{GHASEMIMARZBALI2023120351}& Energy Hub & R2 & - & ECO & - & Implicit & E & - & (MI)NLP (CEN) \\
\citet{GOH2023117484}& Energy Hub & R2 \& R8 & BE & ECO, CUR \& GI  & - & - & E \& HY & E \& HY & Heuristic (CEN) \\
\citet{pedro10161834}& Energy Community & R2 \& R8 & BE & E & - & - & - & E \& G & Deterministic LP (CEN) \\
\citet{GUO202133039}& Industrial Park & R2 & - & ECO & MP & Implicit & E, HE, G \& HY & E, G \& HE & Robust (CEN) \\
\citet{HOU2024130617} & Industrial Park & R2 \& R8 & FR & ECO \& EMI & - & Implicit & E, HE \& HY & E, G \& CA$^M$ & Stackelberg (CEN)  \\
\citet{HWANGBO2022122006}& Industrial Park & R2 & - & ECO & GEN & - & E & - & Stochastic (CEN) \\
\citet{Kantor_Frontier}& Industrial Park & R2, R3 \& R9 & EC \& BE & ECO \& EMI & - & - & - & - & Deterministic MILP (CEN) \\
\citet{KHALIGH2023117354}& Energy Hub & R2 \& R8 & FR & ECO \& EMI & GEN, DEM \& MP & - & E, HY \& C & E, G, HY \& PtX & Stochastic (CEN) \\
\citet{LI2022105261}& Energy Hub & R2, R3 \& R8 & EC \& FR & ENE & GEN \& DEM & - & G \& WAS & G & Robust (CEN) \\
\citet{LI2024119946}& Energy Hub & R2 \& R8 & FR & ECO, EMI \& CUR  & - & Implicit & E, HE \& CO & E \& G & Stackelberg (CEN) \\
\citet{LI2024129784}& Energy Hub& R2 & - & ECO, EMI \& ENE & - & - & E \& HE & E \& G & Heuristic (CEN) \\
\citet{LIN2021117741}& Energy Hub & R2 & - & ECO \& EMI & - & - & E & E \& G & Deterministic MILP (CEN) \\
\citet{LIU2024118204} & Energy Community & R2 & - & ECO & GEN \& DEM & - & - & E & Stochastic (CEN) \\
\citet{LIU2024133528}& Energy Hub & R2, R8 \& R9 & BE & ECO & - & - & E, HE \& G & E, G \& HE & Stochastic (CEN) \\
\citet{MALDET2024100106}& Energy Community & R2 \& R9 & BE & ECO & - & - & HE \& WAT & E & Deterministic MILP (CEN) \\
\citet{MENEGHETTI2012263}& Industrial Park & R2 \& R9 & BE & ECO & - & - & - & E \& G & Deterministic MILP (CEN) \\
\citet{MOHTAVIPOUR2024132880}& Energy Hub & R2 & - & ECO & GEN & Implicit & - & E \& G & Stochastic (CEN) \\
\citet{MOSTAFAVISANI2023116965}& Enery Hub & R2, R8 \& R9 & FR \& BE & ECO \& EMI & MP \& WG & - & E, HE \& HY  & E, HE \& HY & Robust (CEN) \\
\citet{MUGHEES2023121150}& Energy Hub & R2 & - & ECO \& DEG & - & Implicit & E \& HE & E \& HE$^M$ & RL (CEN) \\
\citet{PANG2023125201}& Industrial Park & R2 & - & ECO \& EMI & - & - & E, HE \& HY  & - & (MI)NLP (CEN) \\
\citet{PU_10637267}& Energy Hub & R2 \& R3 & EC & ECO \& EMI & - & Explicit & E \& HY & E$^M$ & Stackelberg (DEC) \\
\citet{su16125016}& Energy Hub & R2 & - & ECO, EMI \& ENE & - & - & E \& HE & E \& G & Heuristic (DEC) \\
\citet{QIAN2024110504}& Indsutrial Park & R2 \& R8 & FR & ECO, EMI \& CUR & - & Explicit & E, HE \& G & E \& G & Multi-layer MILP (CEN) \\
\citet{REN2023119499}& Energy Hub & R2 \& R3 & EC & ECO, EMI \& ENE  & GEN \& DEM & - & HE, G, HY & E, G \& HY$^M$ & Heuristic (CEN) \\
\citet{REN2025135422}& Energy Hub & R2 & - & ECO, EMI \& ENE & - & - & E \& HE & E \& G & Heuristic (DEC) \\
\citet{SALYANI2025106245}& Energy Hub & R2 & - & ECO \& EMI & GEN, DEM \& MP & Explicit & E, HE, G \& HY & E, G \& HY & Robust (CEN) \\
\citet{SOHRABITABAR2022134898}& Energy Hub & R2 & - & ECO & GEN \& MP & Implicit & E, HE, CO \& G & E & Stochastic (CEN) \\
\citet{TOSTADOVELIZ2024123389}& Industrial Park & R2 & - & ECO & GEN, DEM \& MP & - & E & E \& HY & Robust (CEN) \\
\citet{WANG2023126288}& Energy Hub & R2 & - & ECO & GEN & - & - & E \& G$^M$ & Stackelberg (CEN) \\
\citet{en16248041}& Industrial Park & R2 & - & ECO \& EMI & - & - & E \& HE & E & Deterministic MILP (CEN) \\
\citet{WANG20234631}& Energy Hub & R2 & - & ECO, CUR \& DEG & GEN \& MP & - & E \& HY & E \& HY & Stochastic (CEN) \\
\citet{WANG2024129868} & Energy Hub & R2 \& R8 & FR & ECO & - & - & E \& CO & E, G \& HE$^M$ & ADMM (DEC) \\
\citet{WANG2023121141}& Energy Hub & R2 & - & ECO & GEN, DEM \& MP & Explicit & E \& HE & E, G \& HE & Heuristic (CEN) \\
\citet{WEI2022121732}& Industrial Park & R2 & - & ECO & - & - & E \& G & E \& G & Deterministic MILP (CEN) \\
\citet{WEN2025134696}& Energy Hub & R2 \& R3 & EC & ECO \& EMI & GEN \& DEM & Explicit & E, HE \& HY & E, G \& CA & Robust (CEN) \\
\citet{XU2024123525}& Energy Hub & R2 & - & ECO \& EMI & GEN & - & - & E \& HY & Stochastic (CEN) \\
\citet{ZHANG2023127643}& Energy Hub & R2 \& R3 & EC & ECO \& EMI & GEN, DEM \& MP & - & HE, CO \& HY & E \& HY & Chance (CEN) \\
\citet{Zhang_DTU}& Energy Hub & R2, R3 \& R8 & - & ECO \& DEG & GEN \& DEM & - & E \& CO & E & Multi-layer MILP (CEN) \\
\citet{ZHAO2023101340}& Energy Hub & R2, R8 \& R9 & BE & ECO \& EMI & GEN & - & E, HE, G \& WAS & E, G \& HE & Robust (CEN) \\
\citet{ZHENG2025101098}& Energy Community & R2 & - & ECO \& EMI & GEN \& DEM & Implicit & E \& HE & E & Heuristic (CEN) \\
\citet{ZUO2025135413}& Energy Hub & R2 \& R8 & FR & ECO & - & Implicit & HE, CO \& HY & E \& G$^M$ & Stackelberg (CEN) \\
\bottomrule
\end{tabular}
}
\caption{\small Study overview, \textbf{EC}: Energy Cascade, \textbf{FR}: Fuel Replacement, \textbf{BE}: Bio-Energy, \textbf{X}: No operational IS, \textbf{ECO}: Economic, \textbf{EMI}: Emissions, \textbf{CUR}: Curtailment, \textbf{ENE}: Energy, \textbf{DEG}: Degradation, \textbf{GI}: Grid Import, \textbf{DRU}: Demand response utility, \textbf{GEN}: Generation, \textbf{DEM}: Demand, \textbf{MP}: Market Price, \textbf{WQ}: Waste Quality, \textbf{E}: Electricity, \textbf{HE}: Heat, \textbf{CO}: Cooling, \textbf{G}: Gas, \textbf{HY}: Hydrogen, \textbf{WAS}: Waste, \textbf{CA}: Carbon. \textbf{WAT}: Water, \textbf{M}: internal hub market, \textbf{PtX}: Power-to-X Fuel, \textbf{RL}: Reinforcement Learning, \textbf{CEN}: Centralised, \textbf{DEC}: Decentralised}
\label{tab:app}
\end{table}
\end{landscape}

\restoregeometry

%% file: main.bbl
\begin{thebibliography}{92}
\providecommand{\natexlab}[1]{#1}
\providecommand{\url}[1]{\texttt{#1}}
\expandafter\ifx\csname urlstyle\endcsname\relax
  \providecommand{\doi}[1]{doi: #1}\else
  \providecommand{\doi}{doi: \begingroup \urlstyle{rm}\Url}\fi

\bibitem[Abdelaziz et~al.(2011)Abdelaziz, Saidur, and
  Mekhilef]{ABDELAZIZ2011150}
E.A. Abdelaziz, R.~Saidur, and S.~Mekhilef.
\newblock A review on energy saving strategies in industrial sector.
\newblock \emph{Renewable and Sustainable Energy Reviews}, 15\penalty0
  (1):\penalty0 150--168, 2011.
\newblock \doi{https://doi.org/10.1016/j.rser.2010.09.003}.

\bibitem[Alabi et~al.(2022)Alabi, Lu, and Yang]{ALABI2022118997}
Tobi~Michael Alabi, Lin Lu, and Zaiyue Yang.
\newblock Data-driven optimal scheduling of multi-energy system virtual power
  plant (mevpp) incorporating carbon capture system (ccs), electric vehicle
  flexibility, and clean energy marketer (cem) strategy.
\newblock \emph{Applied Energy}, 314:\penalty0 118997, 2022.
\newblock \doi{https://doi.org/10.1016/j.apenergy.2022.118997}.

\bibitem[Albadi and El-Saadany(2008)]{ALBADI20081989}
M.H. Albadi and E.F. El-Saadany.
\newblock A summary of demand response in electricity markets.
\newblock \emph{Electric Power Systems Research}, 78\penalty0 (11):\penalty0
  1989--1996, 2008.
\newblock \doi{https://doi.org/10.1016/j.epsr.2008.04.002}.

\bibitem[ASPIRE(2021)]{P4P_SRIA}
ASPIRE.
\newblock {Processes4Planet - Strategic Research and Innovation Agenda}.
\newblock
  \url{https://www.aspire2050.eu/sites/default/files/users/user85/p4planet_07.06.2022._final.pdf},
  2021.
\newblock [Online; accessed 15-10-2024].

\bibitem[Bocken et~al.(2016)Bocken, de~Pauw, Bakker, and van~der
  Grinten~and]{Bocken03072016}
Nancy M.~P. Bocken, Ingrid de~Pauw, Conny Bakker, and Bram van~der Grinten~and.
\newblock Product design and business model strategies for a circular economy.
\newblock \emph{Journal of Industrial and Production Engineering}, 33\penalty0
  (5):\penalty0 308--320, 2016.
\newblock \doi{https://doi.org/10.1080/21681015.2016.1172124}.

\bibitem[Boyd et~al.(2011)Boyd, Parikh, Chu, Peleato, and
  Eckstein]{boyd2011admm}
Stephen Boyd, Neal Parikh, Eric Chu, Borja Peleato, and Jonathan Eckstein.
\newblock Distributed optimization and statistical learning via the alternating
  direction method of multipliers.
\newblock \emph{Foundations and Trends in Machine Learning}, 3\penalty0
  (1):\penalty0 1--122, 2011.
\newblock \doi{10.1561/2200000016}.

\bibitem[Butturi et~al.(2019)Butturi, Lolli, Sellitto, Balugani, Gamberini, and
  Rimini]{BUTTURI2019113825}
M.A. Butturi, F.~Lolli, M.A. Sellitto, E.~Balugani, R.~Gamberini, and
  B.~Rimini.
\newblock Renewable energy in eco-industrial parks and urban-industrial
  symbiosis: A literature review and a conceptual synthesis.
\newblock \emph{Applied Energy}, 255:\penalty0 113825, 2019.
\newblock \doi{https://doi.org/10.1016/j.apenergy.2019.113825}.

\bibitem[Bänsch et~al.(2021)Bänsch, Busse, Meisel, Rieck, Scholz, Volling,
  and Wichmann]{BANSCH2021107456}
Kristian Bänsch, Jan Busse, Frank Meisel, Julia Rieck, Sebastian Scholz,
  Thomas Volling, and Matthias~G. Wichmann.
\newblock Energy-aware decision support models in production environments: A
  systematic literature review.
\newblock \emph{Computers \& Industrial Engineering}, 159:\penalty0 107456,
  2021.
\newblock \doi{https://doi.org/10.1016/j.cie.2021.107456}.

\bibitem[Cai et~al.(2025)Cai, Wen, Cao, and Qiao]{CAI2025136178}
Pengcheng Cai, Chuanbo Wen, Baosen Cao, and Jinpeng Qiao.
\newblock A wasserstein metric distributionally robust chance-constrained peer
  aggregation energy sharing mechanism for hydrogen-based microgrids
  considering low-carbon drivers.
\newblock \emph{Energy}, 325:\penalty0 136178, 2025.
\newblock \doi{https://doi.org/10.1016/j.energy.2025.136178}.

\bibitem[Chen et~al.(2023{\natexlab{a}})Chen, Sun, Xie, Lin, and
  Wu]{CHEN2023128641}
Minghao Chen, Yi~Sun, Zhiyuan Xie, Nvgui Lin, and Peng Wu.
\newblock An efficient and privacy-preserving algorithm for multiple energy
  hubs scheduling with federated and matching deep reinforcement learning.
\newblock \emph{Energy}, 284:\penalty0 128641, 2023{\natexlab{a}}.
\newblock \doi{https://doi.org/10.1016/j.energy.2023.128641}.

\bibitem[Chen et~al.(2023{\natexlab{b}})Chen, Chang, and Li]{Chen_Wei}
Wei Chen, Xuewu Chang, and Jianing Li.
\newblock A day-ahead optimal scheduling model of an integrated energy system
  for a facility agricultural–industrial park.
\newblock \emph{IET Energy Systems Integration}, 5\penalty0 (3):\penalty0
  261--274, 2023{\natexlab{b}}.
\newblock \doi{https://doi.org/10.1049/esi2.12101}.

\bibitem[Chen et~al.(2024{\natexlab{a}})Chen, Li, Fang, Zhang, Huang, and
  Chen]{chen_iet_nash}
Yanbo Chen, Jiaqi Li, Zhe Fang, Ning Zhang, Tao Huang, and Zuomao Chen.
\newblock Low-carbon optimal scheduling strategy for multi-agent integrated
  energy system of the park based on stackelberg–nash game.
\newblock \emph{IET Energy Systems Integration}, 6\penalty0 (4):\penalty0
  622--637, 2024{\natexlab{a}}.
\newblock \doi{https://doi.org/10.1049/esi2.12164}.

\bibitem[Chen et~al.(2024{\natexlab{b}})Chen, Guo, Du, Yang, and
  Wang]{CHEN2024120734}
Yuzhu Chen, Weimin Guo, Na~Du, Kun Yang, and Jiangjiang Wang.
\newblock Master slave game-based optimization of an off-grid combined cooling
  and power system coupled with solar thermal and photovoltaics considering
  carbon cost allocation.
\newblock \emph{Renewable Energy}, 229:\penalty0 120734, 2024{\natexlab{b}}.
\newblock \doi{https://doi.org/10.1016/j.renene.2024.120734}.

\bibitem[Chinese et~al.(2022)Chinese, Orrù, Meneghetti, Cortella, Giordano,
  and Benedetti]{CHINESE2022124785}
D.~Chinese, P.F. Orrù, A.~Meneghetti, G.~Cortella, L.~Giordano, and
  M.~Benedetti.
\newblock Symbiotic and optimized energy supply for decarbonizing cheese
  production: An italian case study.
\newblock \emph{Energy}, 257:\penalty0 124785, 2022.
\newblock \doi{https://doi.org/10.1016/j.energy.2022.124785}.

\bibitem[Davoudi et~al.(2022)Davoudi, Jooshaki, Moeini-Aghtaie, {Hossein
  Barmayoon}, and Aien]{DAVOUDI2022107889}
Mehdi Davoudi, Mohammad Jooshaki, Moein Moeini-Aghtaie, Mohammad {Hossein
  Barmayoon}, and Morteza Aien.
\newblock Developing a multi-objective multi-layer model for optimal design of
  residential complex energy systems.
\newblock \emph{International Journal of Electrical Power \& Energy Systems},
  138:\penalty0 107889, 2022.
\newblock \doi{https://doi.org/10.1016/j.ijepes.2021.107889}.

\bibitem[Dong and Zhao(2025)]{DONG20251122}
Xiangxiang Dong and Yanling Zhao.
\newblock An economic way to reduce emissions of industrial parks with
  hydrogen-based integrated energy systems.
\newblock \emph{International Journal of Hydrogen Energy}, 106:\penalty0
  1122--1133, 2025.
\newblock \doi{https://doi.org/10.1016/j.ijhydene.2025.01.458}.

\bibitem[Dwijendra et~al.(2022)Dwijendra, Candra, Mubarak, Braiber, Ali, Muda,
  Sivaraman, and Iswanto]{Dwijendra}
Ngakan Ketut~Acwin Dwijendra, Oriza Candra, Ihsan~Ali Mubarak, Hassan~Taher
  Braiber, Muneam~Hussein Ali, Iskandar Muda, R.~Sivaraman, and A.~Heri
  Iswanto.
\newblock Optimal energy-saving in smart energy hub considering demand
  management.
\newblock \emph{Environmental and Climate Technologies}, 26\penalty0
  (1):\penalty0 1244--1256, 2022.
\newblock \doi{doi:10.2478/rtuect-2022-0094}.

\bibitem[{Ellen Macauthur Foundation}(2021)]{EMF_butterfly}
{Ellen Macauthur Foundation}.
\newblock The butterfly diagram: visualising the circular economy.
\newblock
  \url{https://www.ellenmacarthurfoundation.org/circular-economy-diagram},
  2021.
\newblock [Online; accessed 9-12-2024].

\bibitem[Fan et~al.(2024)Fan, Zhang, Wang, Wang, Chen, Ding, and
  Su]{fan10878712}
Wentao Fan, Yongqiang Zhang, Fang Wang, Jianbo Wang, Mengyang Chen, Zeqi Ding,
  and Juan Su.
\newblock Design and optimization of circular economy energy supply and
  utilization model for aquaculture courtyards considering the reuse of manure
  and sewage in organic composting systems.
\newblock In \emph{2024 11th International Forum on Electrical Engineering and
  Automation (IFEEA)}, pages 1145--1151, 2024.
\newblock \doi{https://doi.org/10.1109/ifeea64237.2024.10878712}.

\bibitem[Feigenbaum et~al.(2007)Feigenbaum, Schapira, and
  Shenker]{Feigenbaum2007363}
Joan Feigenbaum, Michael Schapira, and Scott Shenker.
\newblock \emph{Distributed algorithmic mechanism design}, volume
  9780521872829.
\newblock Cambridge University Press, 2007.
\newblock \doi{https://doi.org/10.1017/CBO9780511800481.016}.
\newblock Cited by: 35.

\bibitem[Fraccascia et~al.(2021)Fraccascia, Yazdanpanah, van Capelleveen, and
  Yazan]{Energy_IS}
Luca Fraccascia, Vahid Yazdanpanah, Guido van Capelleveen, and Devrim~Murat
  Yazan.
\newblock Energy-based industrial symbiosis: a literature review for circular
  energy transition.
\newblock \emph{Environment Development and Sustainability}, 23, 04 2021.
\newblock \doi{10.1007/s10668-020-00840-9}.

\bibitem[Freire-Barceló et~al.(2022)Freire-Barceló, Martín-Martínez, and
  Álvaro Sánchez-Miralles]{FREIREBARCELO2022107953}
Teresa Freire-Barceló, Francisco Martín-Martínez, and Álvaro
  Sánchez-Miralles.
\newblock A literature review of explicit demand flexibility providing energy
  services.
\newblock \emph{Electric Power Systems Research}, 209:\penalty0 107953, 2022.
\newblock \doi{https://doi.org/10.1016/j.epsr.2022.107953}.

\bibitem[Galvan-Cara et~al.(2022)Galvan-Cara, Graells, and
  Espuña]{GALVANCARA2022117734}
Aldwin-Lois Galvan-Cara, Moisès Graells, and Antonio Espuña.
\newblock Application of industrial symbiosis principles to the management of
  utility networks.
\newblock \emph{Applied Energy}, 305:\penalty0 117734, 2022.
\newblock \doi{https://doi.org/10.1016/j.apenergy.2021.117734}.

\bibitem[Ghasemi-Marzbali et~al.(2023)Ghasemi-Marzbali, Shafiei, and
  Ahmadiahangar]{GHASEMIMARZBALI2023120351}
Ali Ghasemi-Marzbali, Mohammad Shafiei, and Roya Ahmadiahangar.
\newblock Day-ahead economical planning of multi-vector energy district
  considering demand response program.
\newblock \emph{Applied Energy}, 332:\penalty0 120351, 2023.
\newblock \doi{https://doi.org/10.1016/j.apenergy.2022.120351}.

\bibitem[Gibbs and Deutz(2007)]{GIBBS20071683}
David Gibbs and Pauline Deutz.
\newblock Reflections on implementing industrial ecology through eco-industrial
  park development.
\newblock \emph{Journal of Cleaner Production}, 15\penalty0 (17):\penalty0
  1683--1695, 2007.
\newblock \doi{https://doi.org/10.1016/j.jclepro.2007.02.003}.
\newblock From Material Flow Analysis to Material Flow Management.

\bibitem[Goh et~al.(2023)Goh, Zhang, Ho, and Chew]{GOH2023117484}
Qi~Hao Goh, Lian Zhang, Yong~Kuen Ho, and Irene Mei~Leng Chew.
\newblock Modelling and multi-objective optimisation of sustainable
  solar-biomass-based hydrogen and electricity co-supply hub using
  metaheuristic-topsis approach.
\newblock \emph{Energy Conversion and Management}, 293:\penalty0 117484, 2023.
\newblock \doi{https://doi.org/10.1016/j.enconman.2023.117484}.

\bibitem[Guimarães et~al.(2023)Guimarães, Moreno, Mello, and
  Villar]{pedro10161834}
Pedro Guimarães, Armando Moreno, João Mello, and José Villar.
\newblock A framework for circular energy communities in the agricultural
  sector with a cogeneration case study.
\newblock In \emph{2023 19th International Conference on the European Energy
  Market (EEM)}, pages 1--6, 2023.
\newblock \doi{10.1109/EEM58374.2023.10161834}.

\bibitem[Guo et~al.(2021)Guo, Nojavan, Lei, and Liang]{GUO202133039}
Qun Guo, Sayyad Nojavan, Shi Lei, and Xiaodan Liang.
\newblock Potential evaluation of power-to-hydrogen-to methane conversion
  technology in robust optimal energy management of a multi-energy industrial
  park.
\newblock \emph{International Journal of Hydrogen Energy}, 46\penalty0
  (66):\penalty0 33039--33052, 2021.
\newblock \doi{https://doi.org/10.1016/j.ijhydene.2021.07.148}.

\bibitem[Hong and Gasparatos(2020)]{HONG2020122853}
Hongru Hong and Alexandros Gasparatos.
\newblock Eco-industrial parks in china: Key institutional aspects,
  sustainability impacts, and implementation challenges.
\newblock \emph{Journal of Cleaner Production}, 274:\penalty0 122853, 2020.
\newblock \doi{https://doi.org/10.1016/j.jclepro.2020.122853}.

\bibitem[Hou et~al.(2024)Hou, Ge, Yan, Lu, Zhang, and Dong]{HOU2024130617}
Hui Hou, Xiangdi Ge, Yulin Yan, Yanchao Lu, Ji~Zhang, and Zhao~Yang Dong.
\newblock An integrated energy system “green-carbon” offset mechanism and
  optimization method with stackelberg game.
\newblock \emph{Energy}, 294:\penalty0 130617, 2024.
\newblock \doi{https://doi.org/10.1016/j.energy.2024.130617}.

\bibitem[Hussain et~al.(2024)Hussain, Huang, Li, Zhang, Hussain, Ahmed, and
  Manzoor]{HUSSAIN2024110458}
Jawad Hussain, Qi~Huang, Jian Li, Zhenyuan Zhang, Fazal Hussain, Syed~Adrees
  Ahmed, and Kashif Manzoor.
\newblock Optimization of social welfare in p2p community microgrid with
  efficient decentralized energy management and communication-efficient power
  trading.
\newblock \emph{Journal of Energy Storage}, 81:\penalty0 110458, 2024.
\newblock \doi{https://doi.org/10.1016/j.est.2024.110458}.

\bibitem[Hwangbo et~al.(2022)Hwangbo, Heo, and Yoo]{HWANGBO2022122006}
Soonho Hwangbo, SungKu Heo, and ChangKyoo Yoo.
\newblock Development of deterministic-stochastic model to integrate variable
  renewable energy-driven electricity and large-scale utility networks: Towards
  decarbonization petrochemical industry.
\newblock \emph{Energy}, 238:\penalty0 122006, 2022.
\newblock \doi{https://doi.org/10.1016/j.energy.2021.122006}.

\bibitem[Jeong(2023)]{JEONG2023106190}
In-Jae Jeong.
\newblock A review of decentralized optimization focused on information flows
  of decomposition algorithms.
\newblock \emph{Computers \& Operations Research}, 153:\penalty0 106190, 2023.
\newblock \doi{https://doi.org/10.1016/j.cor.2023.106190}.

\bibitem[Kantor et~al.(2020)Kantor, Robineau, Bütün, and
  Maréchal]{Kantor_Frontier}
Ivan Kantor, Jean-Loup Robineau, Hür Bütün, and François Maréchal.
\newblock A mixed-integer linear programming formulation for optimizing
  multi-scale material and energy integration.
\newblock \emph{Frontiers in Energy Research}, 8, 2020.
\newblock \doi{https://doi.org/10.3389/fenrg.2020.00049}.

\bibitem[Khaligh et~al.(2023)Khaligh, Ghezelbash, Akhtar, Zarei, Liu, and
  Won]{KHALIGH2023117354}
Vahid Khaligh, Azam Ghezelbash, Malik~Sajawal Akhtar, Mohammadamin Zarei, Jay
  Liu, and Wangyun Won.
\newblock Optimal integration of a low-carbon energy system – a circular
  hydrogen economy perspective.
\newblock \emph{Energy Conversion and Management}, 292:\penalty0 117354, 2023.
\newblock \doi{https://doi.org/10.1016/j.enconman.2023.117354}.

\bibitem[Kiani-Moghaddam et~al.(2023)Kiani-Moghaddam, Soltani, Kalogirou,
  Mahian, and Arabkoohsar]{KIANIMOGHADDAM2023128263}
Mohammad Kiani-Moghaddam, Mohsen~N. Soltani, Soteris~A. Kalogirou, Omid Mahian,
  and Ahmad Arabkoohsar.
\newblock A review of neighborhood level multi-carrier energy
  hubs—uncertainty and problem-solving process.
\newblock \emph{Energy}, 281:\penalty0 128263, 2023.
\newblock \doi{https://doi.org/10.1016/j.energy.2023.128263}.

\bibitem[Kirchherr et~al.(2017)Kirchherr, Reike, and Hekkert]{KIRCHHERR2017221}
Julian Kirchherr, Denise Reike, and Marko Hekkert.
\newblock Conceptualizing the circular economy: An analysis of 114 definitions.
\newblock \emph{Resources, Conservation and Recycling}, 127:\penalty0 221--232,
  2017.
\newblock \doi{https://doi.org/10.1016/j.resconrec.2017.09.005}.

\bibitem[Li et~al.(2022)Li, Xu, Zhang, Yuan, and Prokhasko]{LI2022105261}
Ji~Li, Lei Xu, Yuying Zhang, Zhi Yuan, and Lubov~Savelievna Prokhasko.
\newblock Utilization of surplus renewable energy for a gas-fired power plant's
  carbon recycling considering gas storages and imposed regret of
  uncertainties.
\newblock \emph{Journal of Energy Storage}, 54:\penalty0 105261, 2022.
\newblock \doi{https://doi.org/10.1016/j.est.2022.105261}.

\bibitem[Li et~al.(2024{\natexlab{a}})Li, Miao, Lim, Sethanan, and
  Tseng]{LI2024119946}
Ling-Ling Li, Yan Miao, Ming~K. Lim, Kanchana Sethanan, and Ming-Lang Tseng.
\newblock Integrated energy system for low-carbon economic operation
  optimization: Pareto compromise programming and master-slave game.
\newblock \emph{Renewable Energy}, 222:\penalty0 119946, 2024{\natexlab{a}}.
\newblock \doi{https://doi.org/10.1016/j.renene.2024.119946}.

\bibitem[Li et~al.(2024{\natexlab{b}})Li, Qu, Tseng, Lim, Ren, and
  Miao]{LI2024129784}
Ling-Ling Li, Li-Nan Qu, Ming-Lang Tseng, Ming~K. Lim, Xin-Yu Ren, and Yan
  Miao.
\newblock Optimization and performance assessment of solar-assisted combined
  cooling, heating and power system systems: Multi-objective gradient-based
  optimizer.
\newblock \emph{Energy}, 289:\penalty0 129784, 2024{\natexlab{b}}.
\newblock \doi{https://doi.org/10.1016/j.energy.2023.129784}.

\bibitem[Lin et~al.(2021)Lin, Zhong, Wang, Huang, Bai, Wang, Shah, Xie, and
  Zhao]{LIN2021117741}
Jian Lin, Xiaoyi Zhong, Jing Wang, Yuan Huang, Xuetao Bai, Xiaonan Wang, Nilay
  Shah, Shan Xie, and Yingru Zhao.
\newblock Relative optimization potential: A novel perspective to address
  trade-off challenges in urban energy system planning.
\newblock \emph{Applied Energy}, 304:\penalty0 117741, 2021.
\newblock \doi{https://doi.org/10.1016/j.apenergy.2021.117741}.

\bibitem[Liu et~al.(2024{\natexlab{a}})Liu, Xu, Liu, Liu, Hu, Yang, Jawad, and
  Wei]{LIU2024118204}
Yi~Liu, Xiao Xu, Youbo Liu, Junyong Liu, Weihao Hu, Nan Yang, Shafqat Jawad,
  and Zhaobin Wei.
\newblock Scenario-based operation of an integrated rural multi-energy system
  considering agent-based farmer-behavior modeling.
\newblock \emph{Energy Conversion and Management}, 304:\penalty0 118204,
  2024{\natexlab{a}}.
\newblock \doi{https://doi.org/10.1016/j.enconman.2024.118204}.

\bibitem[Liu et~al.(2024{\natexlab{b}})Liu, Xu, Liu, Liu, Hu, Yang, Jawad, and
  Wei]{LIU2024133528}
Yi~Liu, Xiao Xu, Youbo Liu, Junyong Liu, Weihao Hu, Nan Yang, Shafqat Jawad,
  and Zhaobin Wei.
\newblock Operational optimization of a rural multi-energy system supported by
  a joint biomass-solid-waste-energy conversion system and supply chain.
\newblock \emph{Energy}, 312:\penalty0 133528, 2024{\natexlab{b}}.
\newblock \doi{https://doi.org/10.1016/j.energy.2024.133528}.

\bibitem[Maldet et~al.(2024)Maldet, Schwabeneder, Lettner, Loschan,
  Corinaldesi, and Auer]{MALDET2024100106}
Matthias Maldet, Daniel Schwabeneder, Georg Lettner, Christoph Loschan, Carlo
  Corinaldesi, and Hans Auer.
\newblock Local sustainable communities: Sector coupling and community
  optimization in decentralized energy systems.
\newblock \emph{Cleaner Energy Systems}, 7:\penalty0 100106, 2024.
\newblock \doi{https://doi.org/10.1016/j.cles.2023.100106}.

\bibitem[Meneghetti and Nardin(2012)]{MENEGHETTI2012263}
Antonella Meneghetti and Gioacchino Nardin.
\newblock Enabling industrial symbiosis by a facilities management optimization
  approach.
\newblock \emph{Journal of Cleaner Production}, 35:\penalty0 263--273, 2012.
\newblock \doi{https://doi.org/10.1016/j.jclepro.2012.06.002}.

\bibitem[Mohammadi et~al.(2017)Mohammadi, Noorollahi, Mohammadi-ivatloo, and
  Yousefi]{MOHAMMADI20171512}
Mohammad Mohammadi, Younes Noorollahi, Behnam Mohammadi-ivatloo, and Hossein
  Yousefi.
\newblock Energy hub: From a model to a concept – a review.
\newblock \emph{Renewable and Sustainable Energy Reviews}, 80:\penalty0
  1512--1527, 2017.
\newblock \doi{https://doi.org/10.1016/j.rser.2017.07.030}.

\bibitem[Mohtavipour(2024)]{MOHTAVIPOUR2024132880}
Seyed~Saeid Mohtavipour.
\newblock Convex relaxation of two-stage network-constrained stochastic
  programming for chp microgrid optimal scheduling.
\newblock \emph{Energy}, 308:\penalty0 132880, 2024.
\newblock \doi{https://doi.org/10.1016/j.energy.2024.132880}.

\bibitem[{Mostafavi Sani} et~al.(2023){Mostafavi Sani}, Afshari, and
  Saif]{MOSTAFAVISANI2023116965}
Mostafa {Mostafavi Sani}, Hamid Afshari, and Ahmed Saif.
\newblock A robust framework for waste-to-energy technology selection: A case
  study in nova scotia, canada.
\newblock \emph{Energy Conversion and Management}, 284:\penalty0 116965, 2023.
\newblock \doi{https://doi.org/10.1016/j.enconman.2023.116965}.

\bibitem[Mughees et~al.(2023)Mughees, Jaffery, Mughees, Ansari, and
  Mughees]{MUGHEES2023121150}
Neelam Mughees, Mujtaba~Hussain Jaffery, Anam Mughees, Ejaz~Ahmad Ansari, and
  Abdullah Mughees.
\newblock Reinforcement learning-based composite differential evolution for
  integrated demand response scheme in industrial microgrids.
\newblock \emph{Applied Energy}, 342:\penalty0 121150, 2023.
\newblock \doi{https://doi.org/10.1016/j.apenergy.2023.121150}.

\bibitem[M{\"u}nster et~al.(2020)M{\"u}nster, Sneum, Pedersen, B{\"u}hler,
  Elmegaard, Giannelos, Zhang, Strbac, Berger, Radu, Elsaesser, Oudalov, and
  Iliceto]{White_sector}
Marie M{\"u}nster, {Daniel M{\o}ller} Sneum, {Rasmus Bramstoft} Pedersen,
  Fabian B{\"u}hler, Brian Elmegaard, Spyros Giannelos, Xi~Zhang, Goran Strbac,
  Mathias Berger, David Radu, Damian Elsaesser, Alexandre Oudalov, and Antonio
  Iliceto.
\newblock \emph{Sector Coupling: Concepts, State-of-the-art and Perspectives}.
\newblock European Technology and Innovation Platform, 2020.

\bibitem[Nozari et~al.(2022)Nozari, Yaghoubi, Jafarpur, and
  Mansoori]{NOZARI2022103972}
Mohammad~Hossein Nozari, Mahmoud Yaghoubi, Khosrow Jafarpur, and G.~Ali
  Mansoori.
\newblock Development of dynamic energy storage hub concept: A comprehensive
  literature review of multi storage systems.
\newblock \emph{Journal of Energy Storage}, 48:\penalty0 103972, 2022.
\newblock \doi{https://doi.org/10.1016/j.est.2022.103972}.

\bibitem[Nozarian et~al.(2024)Nozarian, Hajizadeh, and
  Fereidunian]{NOZARIAN2024103834}
Mahdi Nozarian, Amin Hajizadeh, and Alireza Fereidunian.
\newblock A methodological review on management of building cluster meso energy
  hubs in accordance with the agile framework: Exploring flexibility upward the
  building-cluster-city hierarchy.
\newblock \emph{Sustainable Energy Technologies and Assessments}, 67:\penalty0
  103834, 2024.
\newblock \doi{https://doi.org/10.1016/j.seta.2024.103834}.

\bibitem[Page et~al.(2021)Page, Moher, Bossuyt, Boutron, Hoffmann, Mulrow,
  Shamseer, Tetzlaff, Akl, Brennan, Chou, Glanville, Grimshaw,
  Hr{\'o}bjartsson, Lalu, Li, Loder, Mayo-Wilson, McDonald, McGuinness,
  Stewart, Thomas, Tricco, Welch, Whiting, and McKenzie]{Pagen160}
Matthew~J Page, David Moher, Patrick~M Bossuyt, Isabelle Boutron, Tammy~C
  Hoffmann, Cynthia~D Mulrow, Larissa Shamseer, Jennifer~M Tetzlaff, Elie~A
  Akl, Sue~E Brennan, Roger Chou, Julie Glanville, Jeremy~M Grimshaw,
  Asbj{\o}rn Hr{\'o}bjartsson, Manoj~M Lalu, Tianjing Li, Elizabeth~W Loder,
  Evan Mayo-Wilson, Steve McDonald, Luke~A McGuinness, Lesley~A Stewart, James
  Thomas, Andrea~C Tricco, Vivian~A Welch, Penny Whiting, and Joanne~E
  McKenzie.
\newblock Prisma 2020 explanation and elaboration: updated guidance and
  exemplars for reporting systematic reviews.
\newblock \emph{BMJ}, 372, 2021.
\newblock \doi{https://doi.org/10.1136/bmj.n160}.

\bibitem[Pang et~al.(2023)Pang, Liew, Woon, Ho, {Wan Alwi}, and
  Klemeš]{PANG2023125201}
Kang~Ying Pang, Peng~Yen Liew, Kok~Sin Woon, Wai~Shin Ho, Sharifah~Rafidah {Wan
  Alwi}, and Jiří~Jaromír Klemeš.
\newblock Multi-period multi-objective optimisation model for multi-energy
  urban-industrial symbiosis with heat, cooling, power and hydrogen demands.
\newblock \emph{Energy}, 262:\penalty0 125201, 2023.
\newblock \doi{https://doi.org/10.1016/j.energy.2022.125201}.

\bibitem[Patala et~al.(2020)Patala, Salmi, and Bocken]{PATALA2020121093}
Samuli Patala, Asta Salmi, and Nancy Bocken.
\newblock Intermediation dilemmas in facilitated industrial symbiosis.
\newblock \emph{Journal of Cleaner Production}, 261:\penalty0 121093, 2020.
\newblock \doi{https://doi.org/10.1016/j.jclepro.2020.121093}.

\bibitem[Pu et~al.(2025)Pu, Liu, Wang, and Liu]{PU_10637267}
Yue Pu, Haoming Liu, Jian Wang, and Feng Liu.
\newblock Decentralized coordinated scheduling of port integrated energy system
  and bulk terminal with two-tier incentives.
\newblock \emph{IEEE Transactions on Smart Grid}, 16\penalty0 (1):\penalty0
  101--114, 2025.
\newblock \doi{10.1109/TSG.2024.3444276}.

\bibitem[Qi et~al.(2024)Qi, Peng, Wu, and Tseng]{su16125016}
Zhaoyu Qi, Shitao Peng, Peisen Wu, and Ming-Lang Tseng.
\newblock Renewable energy distributed energy system optimal configuration and
  performance analysis: Improved zebra optimization algorithm.
\newblock \emph{Sustainability}, 16\penalty0 (12), 2024.
\newblock \doi{https://doi.org/10.3390/su16125016}.

\bibitem[Qian et~al.(2024)Qian, Lin, Li, Wang, and Li]{QIAN2024110504}
Liang Qian, Shunfu Lin, Fangxing Li, Wei Wang, and Dongdong Li.
\newblock Low carbon optimization dispatching of energy intensive industrial
  park based on adaptive stepped demand response incentive mechanism.
\newblock \emph{Electric Power Systems Research}, 233:\penalty0 110504, 2024.
\newblock \doi{https://doi.org/10.1016/j.epsr.2024.110504}.

\bibitem[{Ramir D.T. Certeza} et~al.(2025){Ramir D.T. Certeza}, Purnama, Ahsan,
  Low, and Lu]{RAMIRDTCERTEZA2025115377}
La~Verne {Ramir D.T. Certeza}, Aloisius~Rabata Purnama, Aniq Ahsan,
  Jonathan~S.C. Low, and Wen~F. Lu.
\newblock Review of mathematical programming models for energy-based industrial
  symbiosis networks.
\newblock \emph{Renewable and Sustainable Energy Reviews}, 212:\penalty0
  115377, 2025.
\newblock \doi{https://doi.org/10.1016/j.rser.2025.115377}.

\bibitem[Ren et~al.(2023)Ren, Li, and Li]{REN2023119499}
Ting Ren, Ran Li, and Xin Li.
\newblock Bi-level multi-objective robust optimization for performance
  improvements in integrated energy system with solar fuel production.
\newblock \emph{Renewable Energy}, 219:\penalty0 119499, 2023.
\newblock \doi{https://doi.org/10.1016/j.renene.2023.119499}.

\bibitem[Ren et~al.(2025)Ren, Lei, Li, Guo, Chen, Yap, and Wang]{REN2025135422}
Xinyu Ren, Haowen Lei, Yixuan Li, Xiaolong Guo, Zhonghao Chen, Pow-Seng Yap,
  and Zhihua Wang.
\newblock A multi-criteria assessment method for design and dispatch of
  distributed energy systems considering different energy consumption
  attributes.
\newblock \emph{Energy}, 323:\penalty0 135422, 2025.
\newblock \doi{https://doi.org/10.1016/j.energy.2025.135422}.

\bibitem[Salyani et~al.(2025)Salyani, Zare, Javani, and
  Boynuegri]{SALYANI2025106245}
Pouya Salyani, Kazem Zare, Nader Javani, and Ali~Rifat Boynuegri.
\newblock Risk-based scheduling of multi-energy microgrids with power-to-x
  technology under a multi-objective framework.
\newblock \emph{Sustainable Cities and Society}, 122:\penalty0 106245, 2025.
\newblock \doi{https://doi.org/10.1016/j.scs.2025.106245}.

\bibitem[Shabanpour-Haghighi and Karimaghaei(2022)]{SHABANPOURHAGHIGHI20226164}
Amin Shabanpour-Haghighi and Mina Karimaghaei.
\newblock An overview on multi-carrier energy networks: From a concept to
  future trends and challenges.
\newblock \emph{International Journal of Hydrogen Energy}, 47\penalty0
  (9):\penalty0 6164--6186, 2022.
\newblock \doi{https://doi.org/10.1016/j.ijhydene.2021.11.257}.

\bibitem[Sharma et~al.(2022)Sharma, {Dutt Mathur}, Mishra, and
  Bansal]{SHARMA2022120028}
Pavitra Sharma, Hitesh {Dutt Mathur}, Puneet Mishra, and Ramesh~C. Bansal.
\newblock A critical and comparative review of energy management strategies for
  microgrids.
\newblock \emph{Applied Energy}, 327:\penalty0 120028, 2022.
\newblock \doi{https://doi.org/10.1016/j.apenergy.2022.120028}.

\bibitem[Sheu and Chen(2014)]{SHEU201465}
Jiuh-Biing Sheu and Yenming~J. Chen.
\newblock Transportation and economies of scale in recycling low-value
  materials.
\newblock \emph{Transportation Research Part B: Methodological}, 65:\penalty0
  65--76, 2014.
\newblock \doi{https://doi.org/10.1016/j.trb.2014.03.008}.

\bibitem[{Sohrabi Tabar} et~al.(2022){Sohrabi Tabar}, Ghassemzadeh, and
  Tohidi]{SOHRABITABAR2022134898}
Vahid {Sohrabi Tabar}, Saeid Ghassemzadeh, and Sajjad Tohidi.
\newblock Risk-based day-ahead planning of a renewable multi-carrier system
  integrated with multi-level electric vehicle charging station, cryptocurrency
  mining farm and flexible loads.
\newblock \emph{Journal of Cleaner Production}, 380:\penalty0 134898, 2022.
\newblock \doi{https://doi.org/10.1016/j.jclepro.2022.134898}.

\bibitem[Song et~al.(2018)Song, Geng, Dong, and Chen]{SONG2018414}
Xiaoqian Song, Yong Geng, Huijuan Dong, and Wei Chen.
\newblock Social network analysis on industrial symbiosis: A case of gujiao
  eco-industrial park.
\newblock \emph{Journal of Cleaner Production}, 193:\penalty0 414--423, 2018.
\newblock \doi{https://doi.org/10.1016/j.jclepro.2018.05.058}.

\bibitem[Tas et~al.(2021)Tas, Tyrkko, Susan, Demir~Duru, Kechichian, Boll, and
  Felix]{EIP_WBG}
Nilgun Tas, Klaus Tyrkko, Christian~R. Susan, Sinem Demir~Duru, Etienne~Raffi
  Kechichian, Mareike Boll, and Steffen Felix.
\newblock An international framework for eco-industrial parks : Version 2.0
  (english).
\newblock \emph{World Bank Group}, 2021.
\newblock URL
  \url{http://documents.worldbank.org/curated/en/350851612561873572/An-International-Framework-for-Eco-Industrial-Parks-Version-2-0}.
\newblock [Online; accessed 15-10-2024].

\bibitem[Tleuken et~al.(2025)Tleuken, Rogetzer, Fraccascia, and
  Yazan]{TLEUKEN2025125324}
Aidana Tleuken, Patricia Rogetzer, Luca Fraccascia, and Devrim~Murat Yazan.
\newblock Designing a stakeholder engagement framework with critical success
  factors for hubs for circularity.
\newblock \emph{Journal of Environmental Management}, 384:\penalty0 125324,
  2025.
\newblock \doi{https://doi.org/10.1016/j.jenvman.2025.125324}.

\bibitem[Tostado-Véliz et~al.(2024)Tostado-Véliz, {Rezaee Jordehi}, Mansouri,
  Escámez, Alharthi, and Jurado]{TOSTADOVELIZ2024123389}
Marcos Tostado-Véliz, Ahmad {Rezaee Jordehi}, Seyed~Amir Mansouri, Antonio
  Escámez, Yahya~Z. Alharthi, and Francisco Jurado.
\newblock Risk-averse electrolyser sizing in industrial parks: An efficient
  stochastic-robust approach.
\newblock \emph{Applied Energy}, 367:\penalty0 123389, 2024.
\newblock \doi{https://doi.org/10.1016/j.apenergy.2024.123389}.

\bibitem[Turken and Geda(2020)]{TURKEN2020104974}
Nazli Turken and Avinash Geda.
\newblock Supply chain implications of industrial symbiosis: A review and
  avenues for future research.
\newblock \emph{Resources, Conservation and Recycling}, 161:\penalty0 104974,
  2020.
\newblock \doi{https://doi.org/10.1016/j.resconrec.2020.104974}.

\bibitem[{United Nations Environment Programme }(2024)]{UNEP}
{United Nations Environment Programme }.
\newblock Global resourse outlook 2024 - bend the trend pathways to a liveable
  planet as resource use spikes.
\newblock
  \url{https://www.resourcepanel.org/sites/default/files/documents/document/media/gro24_full_report_29feb_final_for_web.pdf},
  2024.
\newblock [wedocs.unep.org/20.500.11822/44901, Online; accessed 9-12-2024].

\bibitem[Wang et~al.(2023{\natexlab{a}})Wang, Sayed, Zhang, Zhang, Ren, Jia,
  and Bi]{WANG2023126288}
Cheng Wang, Ahmed~Rabee Sayed, Han Zhang, Xian Zhang, Jianpeng Ren, Qiyue Jia,
  and Tianshu Bi.
\newblock Two-stage distributionally robust strategic offering in pool-based
  coupled electricity and gas market.
\newblock \emph{Energy}, 265:\penalty0 126288, 2023{\natexlab{a}}.
\newblock \doi{https://doi.org/10.1016/j.energy.2022.126288}.

\bibitem[Wang et~al.(2023{\natexlab{b}})Wang, Lu, Zhang, Yan, and
  Feng]{en16248041}
Jiaying Wang, Chunguang Lu, Shuai Zhang, Huajiang Yan, and Changsen Feng.
\newblock Optimal energy management strategy of clustered industry factories
  considering carbon trading and supply chain coupling.
\newblock \emph{Energies}, 16\penalty0 (24), 2023{\natexlab{b}}.
\newblock \doi{https://doi.org/10.3390/en16248041}.

\bibitem[Wang et~al.(2023{\natexlab{c}})Wang, Xu, Wu, Huang, Chen, and
  Hu]{WANG20234631}
Jun Wang, Xiao Xu, Lan Wu, Qi~Huang, Zhe Chen, and Weihao Hu.
\newblock Risk-averse based optimal operational strategy of grid-connected
  photovoltaic/wind/battery/diesel hybrid energy system in the
  electricity/hydrogen markets.
\newblock \emph{International Journal of Hydrogen Energy}, 48\penalty0
  (12):\penalty0 4631--4648, 2023{\natexlab{c}}.
\newblock \doi{https://doi.org/10.1016/j.ijhydene.2022.11.006}.

\bibitem[Wang et~al.(2024)Wang, Xian, Jiao, Liu, Xing, and
  Wang]{WANG2024129868}
L.L. Wang, R.C. Xian, P.H. Jiao, X.H. Liu, Y.W. Xing, and W.~Wang.
\newblock Cooperative operation of industrial/commercial/residential integrated
  energy system with hydrogen energy based on nash bargaining theory.
\newblock \emph{Energy}, 288:\penalty0 129868, 2024.
\newblock \doi{https://doi.org/10.1016/j.energy.2023.129868}.

\bibitem[Wang et~al.(2023{\natexlab{d}})Wang, Lv, and Zeman]{WANG2023121141}
Qiu-Yu Wang, Xian-Long Lv, and Abdol Zeman.
\newblock Optimization of a multi-energy microgrid in the presence of energy
  storage and conversion devices by using an improved gray wolf algorithm.
\newblock \emph{Applied Thermal Engineering}, 234:\penalty0 121141,
  2023{\natexlab{d}}.
\newblock \doi{https://doi.org/10.1016/j.applthermaleng.2023.121141}.

\bibitem[Wang et~al.(2020)Wang, Liu, Liu, and Liu]{WANG2020106060}
Xiaodi Wang, Youbo Liu, Chang Liu, and Junyong Liu.
\newblock Coordinating energy management for multiple energy hubs: From a
  transaction perspective.
\newblock \emph{International Journal of Electrical Power \& Energy Systems},
  121:\penalty0 106060, 2020.
\newblock \doi{https://doi.org/10.1016/j.ijepes.2020.106060}.

\bibitem[Wei et~al.(2022)Wei, Qiu, Liang, Liao, Klemeš, Xue, and
  Zhang]{WEI2022121732}
Xintong Wei, Rui Qiu, Yongtu Liang, Qi~Liao, Jiří~Jaromír Klemeš, Jinjun
  Xue, and Haoran Zhang.
\newblock Roadmap to carbon emissions neutral industrial parks: Energy,
  economic and environmental analysis.
\newblock \emph{Energy}, 238:\penalty0 121732, 2022.
\newblock \doi{https://doi.org/10.1016/j.energy.2021.121732}.

\bibitem[Wen et~al.(2025)Wen, Jia, Cao, Guo, Jiao, Li, and Li]{WEN2025134696}
Jiaxing Wen, Rong Jia, Ge~Cao, Yi~Guo, Yang Jiao, Wei Li, and Peihang Li.
\newblock Robust economic scheduling model for virtual power plant considering
  electrolysis of molten carbonate and dynamic compensation mechanism.
\newblock \emph{Energy}, 317:\penalty0 134696, 2025.
\newblock \doi{https://doi.org/10.1016/j.energy.2025.134696}.

\bibitem[{World Economic Forum} et~al.(2025){World Economic Forum}, {University
  of Cambridge Institute for Sustainability Leadership}, and {Bain \&
  Company}]{wef2025circular}
{World Economic Forum}, {University of Cambridge Institute for Sustainability
  Leadership}, and {Bain \& Company}.
\newblock Circular transformation of industries: Unlocking new value in a
  resource-constrained world.
\newblock Technical report, World Economic Forum, 2025.
\newblock URL
  \url{https://reports.weforum.org/docs/WEF_Circular_Transformation_of_Industries_2025.pdf}.
\newblock Accessed on May 27, 2025.

\bibitem[Xavier et~al.(2024)Xavier, Thollander, Hilletofth, and
  Johansson]{Xavier_SCM_EM}
Bruna~Maria Xavier, Patrik Thollander, Per Hilletofth, and Maria Johansson.
\newblock Exploring energy management integration into upstream supply chains:
  a systematic literature review.
\newblock \emph{Frontiers in Energy Research}, 12, 2024.
\newblock \doi{10.3389/fenrg.2024.1425795}.

\bibitem[Xu et~al.(2024)Xu, Long, Zhao, and Du]{XU2024123525}
Tiantian Xu, Jian Long, Liang Zhao, and Wenli Du.
\newblock Material and energy coupling systems optimization for large-scale
  industrial refinery with sustainable energy penetration under multiple
  uncertainties using two-stage stochastic programming.
\newblock \emph{Applied Energy}, 371:\penalty0 123525, 2024.
\newblock \doi{https://doi.org/10.1016/j.apenergy.2024.123525}.

\bibitem[Yazici et~al.(2022)Yazici, Alakaş, and Eren]{yazici_OR_IS}
Emre Yazici, Hacı Alakaş, and Tamer Eren.
\newblock Analysis of operations research methods for decision problems
  solution in the industrial symbiosis: A literature review.
\newblock \emph{Environmental Science and Pollution Research}, 04 2022.
\newblock \doi{10.21203/rs.3.rs-1590756/v1}.

\bibitem[{Zarei Golambahri} et~al.(2024){Zarei Golambahri}, Shakarami, and
  Doostizadeh]{ZAREIGOLAMBAHRI2024109901}
Milad {Zarei Golambahri}, Mahmoudreza Shakarami, and Meysam Doostizadeh.
\newblock Security-aware joint energy and flexibility trading in
  electricity-heat networks: A novel clearing and validation analysis.
\newblock \emph{International Journal of Electrical Power \& Energy Systems},
  157:\penalty0 109901, 2024.
\newblock \doi{https://doi.org/10.1016/j.ijepes.2024.109901}.

\bibitem[Zhang et~al.(2023)Zhang, Wang, Zhao, Yang, and {Bu
  sinnah}]{ZHANG2023127643}
Hui Zhang, Jiye Wang, Xiongwen Zhao, Jingqi Yang, and Zainab~Ali {Bu sinnah}.
\newblock Modeling a hydrogen-based sustainable multi-carrier energy system
  using a multi-objective optimization considering embedded joint chance
  constraints.
\newblock \emph{Energy}, 278:\penalty0 127643, 2023.
\newblock \doi{https://doi.org/10.1016/j.energy.2023.127643}.

\bibitem[Zhang et~al.(2021)Zhang, Zhou, Wu, Cao, Liu, Voropai, and
  Barakhtenko]{Zhang_DTU}
Kuan Zhang, Bin Zhou, Qiuwei Wu, Yijia Cao, Nian Liu, Nikolai Voropai, and
  Evgeny Barakhtenko.
\newblock Modeling and utilization of biomass-to-syngas for industrial
  multi-energy systems.
\newblock \emph{CSEE Journal of Power and Energy Systems}, 7\penalty0
  (5):\penalty0 932--942, 2021.
\newblock \doi{https://doi.org/10.17775/CSEEJPES.2020.06190}.

\bibitem[Zhao et~al.(2023)Zhao, Cao, Wang, Li, and Wang]{ZHAO2023101340}
Huiru Zhao, Yiqiong Cao, Xuejie Wang, Bingkang Li, and Yuwei Wang.
\newblock Distributionally robust optimization for dispatching the integrated
  operation of rural multi-energy system and waste disposal facilities under
  uncertainties.
\newblock \emph{Energy for Sustainable Development}, 77:\penalty0 101340, 2023.
\newblock \doi{https://doi.org/10.1016/j.esd.2023.101340}.

\bibitem[Zheng and Wu(2025)]{ZHENG2025101098}
Hong Zheng and Zhixin Wu.
\newblock Enhancing economic and environmental performance of energy
  communities: A multi-objective optimization approach with mountain gazelle
  optimizer.
\newblock \emph{Sustainable Computing: Informatics and Systems}, 46:\penalty0
  101098, 2025.
\newblock \doi{https://doi.org/10.1016/j.suscom.2025.101098}.

\bibitem[Zhou et~al.(2018)Zhou, Wu, and Long]{ZHOU2018993}
Yue Zhou, Jianzhong Wu, and Chao Long.
\newblock Evaluation of peer-to-peer energy sharing mechanisms based on a
  multiagent simulation framework.
\newblock \emph{Applied Energy}, 222:\penalty0 993--1022, 2018.
\newblock \doi{https://doi.org/10.1016/j.apenergy.2018.02.089}.

\bibitem[Zhu et~al.(2022)Zhu, Yang, Liu, Wang, Ma, and Guan]{ZHU2022118636}
Dafeng Zhu, Bo~Yang, Yuxiang Liu, Zhaojian Wang, Kai Ma, and Xinping Guan.
\newblock Energy management based on multi-agent deep reinforcement learning
  for a multi-energy industrial park.
\newblock \emph{Applied Energy}, 311:\penalty0 118636, 2022.
\newblock \doi{https://doi.org/10.1016/j.apenergy.2022.118636}.

\bibitem[Zuo et~al.(2025)Zuo, Xi, and Zhang]{ZUO2025135413}
Lujie Zuo, Yufei Xi, and Jiansheng Zhang.
\newblock Leveraging electrochemical {CO2} reduction for optimizing
  comprehensive benefits of multi-energy systems: A collaborative optimization
  approach driven by energy-carbon integrated pricing.
\newblock \emph{Energy}, 322:\penalty0 135413, 2025.
\newblock \doi{https://doi.org/10.1016/j.energy.2025.135413}.

\end{thebibliography}
